\documentclass{article}
\usepackage{mathrsfs}
\usepackage{latexsym}
\usepackage{graphicx}
\usepackage{amssymb}
\usepackage{amsmath}
\usepackage{calrsfs}

\begin{document}
\title{\normalsize {ON THE REGULARITY OF THE SOLUTIONS FOR CAUCHY PROBLEM \\  \vskip0.2cm
OF INCOMPRESSIBLE 3D NAVIER-STOKES EQUATION}}

\author {\normalsize{Qun Lin}  \footnotemark \thanks{School of Mathematical Sciences, Xiamen University, P.
R. China ( e-mail: \qquad \qquad \qquad linqun@xmu.edu.cn ) } }

\date{\begin{footnotesize}September 8, 2012 \end{footnotesize}}
\maketitle
\bigskip

\begin{minipage}{11cm}\begin{footnotesize} {\bf Abstract.} In this paper we will prove that the vorticity
belongs to $L^\infty (0,T;\;L^2({\mathbb R}^3))$ for the Cauchy
problem of 3D incompressible Navier-Stokes equation, then the
existence of a global smooth solution is obtained. Our approach is
to construct a set of auxiliary problems to approximate the original
one of vorticity equation.


\end{footnotesize}

\end{minipage}

\bigskip

\textbf{1. Introduction}

Let $\mathscr{D} ({\mathbb R}^3)$ be the space of $C^\infty $
functions with compact support contained in ${\mathbb R}^3$. Some
basic spaces will be used in this paper:
\begin{equation*}
\begin{split}
 &{\cal V}=\{\,u\in \mathscr{D} ({\mathbb R}^3),\;\;\mbox{div}u=0\,\} \\
 &V=\mbox{the}\;\mbox{closure}\;\mbox{of}\;{\cal
 V}\;\mbox{in}\;H^1({\mathbb
R}^3) \\
 &H=\mbox{the}\;\mbox{closure}\;\mbox{of}\;{\cal
 V}\;\mbox{in}\;L^2({\mathbb
R}^3) \\
\end{split}
\end{equation*}

The velocity-pressure form for Navier- Stokes equation is
\begin{equation}
\begin{split}
 &\partial _t u_1 + u_1 \partial _{x_1 } u_1 + u_2 \partial _{x_2 } u_1
+ u_3 \partial _{x_3 } u_1 + \partial _{x_1 } p=\Delta u_1 \\
 &\partial _t u_2 + u_1 \partial _{x_1 } u_2 + u_2 \partial _{x_2 } u_2 + u_3
\partial _{x_3 } u_2 + \partial _{x_2 } p=\Delta u_2 \\
 &\partial _t u_3 + u_1 \partial _{x_1 } u_3 + u_2 \partial _{x_2 } u_3 + u_3
\partial _{x_3 } u_3 + \partial _{x_3 } p=\Delta u_3 \\
\end{split}
\end{equation}
with the initial conditions $\left. {(u_1 ,u_2 ,u_3 )} \right|_{t=0}
=(u_{10} ,u_{20} ,u_{30} )(x)$, henceforth we always ignore the assumption of sufficient smoothness of the initial conditions.
Moreover,  the incompressible condition is
\[
\partial _{x_1 } u_1 +\partial _{x_2 } u_2 +\partial _{x_3 } u_3 =0
\]
where $x=(x_1,x_2,x_3)$ is a point of ${\mathbb R}^3$,
$u=(u_1,u_2,u_3)$ is velocity, $p$ is pressure, and $\nu>0$ is
viscosity.
\\

The
vorticity-velocity form for Navier-Stokes equation is
\begin{equation}
\begin{split}
 &\partial _t \omega _1 +u_1 \partial _{x_1 } \omega _1 +u_2 \partial
_{x_2 } \omega _1 +u_3 \partial _{x_3 } \omega _1 -\omega _1
\partial _{x_1 } u_1 -\omega _2 \partial _{x_2 } u_1 -\omega _3
\partial _{x_3 } u_1
=\Delta \omega _1 \\
 &\partial _t \omega _2 +u_1 \partial _{x_1 } \omega _2 +u_2 \partial _{x_2 }
\omega _2 +u_3 \partial _{x_3 } \omega _2 -\omega _1 \partial _{x_1 } u_2
-\omega _2 \partial _{x_2 } u_2 -\omega _3 \partial _{x_3 } u_2 =\Delta
\omega _2 \\
 &\partial _t \omega _3 +u_1 \partial _{x_1 } \omega _3 +u_2 \partial _{x_2 }
\omega _3 +u_3 \partial _{x_3 } \omega _3 -\omega _1 \partial _{x_1
} u_3 -\omega _2 \partial _{x_2 } u_3 -\omega _3 \partial _{x_3 }
u_3 =\Delta
\omega _3 \\
\end{split}
\end{equation}
with the initial conditions $\left. {(\omega _1 ,\omega _2 ,\omega _3 )}
\right|_{t=0} =(\omega _{10} ,\omega _{20} ,\omega _{30}
)=(\mbox{curl}u_{10} ,\;\mbox{curl}u_{20} ,\;\mbox{curl}u_{30} )$, and the
incompressible condition :
\begin{equation*}
\begin{split}
 &\partial _{x_1 } u_1 \;+\partial _{x_2 } u_2 \;\,+\partial _{x_3 } u_3 =0
\\
&\partial _{x_1 } \omega _1 +\partial _{x_2 } \omega _2 +\partial
_{x_3 }
\omega _3 =0 \\
\end{split}
\end{equation*}
\\

We here recall the global $L^2$-estimate from [4] for the Navier-Stokes equation of velocity-pressure form.

In the sequel, it is assumed that the initial value $u_0 $ satisfies
the following conditions:
\begin{equation}
\begin{split}
\left| {\partial _{x_j }^\mu u_{i0} (x)} \right| \le C_\mu (1 +
\left| x \right|)^{ - \sigma }, \qquad i,j=1,2,3
\end{split}
\end{equation}
where $ \mu = 0,1 $ and $ \sigma >0 $ is integer.

For the handling the initial value problem, a weighted function is
introduced:
\begin{equation*}
\begin{split}
\theta_{r} =\left\{ {{\begin{array}{*{20}c}
 &{e^{-\;\frac{\left| x \right|^2 }{r^2-\left| x \right|^2}}\quad \left| x \right|<r}
\hfill \\
 &{\quad \;0 \qquad \quad\; \left| x \right|\ge r} \hfill \\
\end{array} }} \right.\quad \quad (r>0)
\end{split}
\end{equation*}
which is of the properties:
\begin{equation}
\begin{split}
\theta _{r} \to 1,\quad \quad \partial _i \theta _{r} \to 0,\quad
\quad \partial _i \partial _j \theta _{r} \to 0
\end{split}
\end{equation}
as $ r\to +\infty $ for each relatively fixed $ x \in \mathbb R^3 $. \\

Moreover, let $v=\theta _{r} u$, we still have
\begin{equation}
\begin{split}
 &\partial _i v=u\,\partial _i \theta _{r} +\theta _{r}
\,\partial _i u \\
 &\partial _i^2 v=u\,\partial _i^2 \theta _{r} +2\,\partial _i
\theta _{r} \partial _i u+\theta _{r} \,\partial _i^2 u \\
 &\partial _i \partial _j v=u\,\partial _i \partial _j \theta _{r}
+\partial _j \theta _{r} \partial _i u+\partial _i \theta _{r} \partial _j u+\theta _{r} \,\partial _i \partial _j u \\
\end{split}
\end{equation}
Since
\begin{equation*}
\begin{split}
 &\int_{{\mathbb R}^3} {\theta _{r} u_i (u_1 \partial _{x_1 } u_i +u_2
\partial _{x_2 } u_i +u_3 \partial _{x_3 } u_i )} =\frac{1}{2}\int_{{\mathbb
R}^3} {\theta _{r} (u_1 \partial _{x_1 } u_i^2 +u_2
\partial _{x_2
} u_i^2 +u_3 \partial _{x_3 } u_i^2 )} \\
 &=-\frac{1}{2}\int_{{\mathbb R}^3} {u_i^2 (\partial _{x_1 } (\theta _{r} u_1 )+\partial _{x_2 } (\theta _{r} u_2 )+\partial
_{x_3 }
(\theta _{r} u_3 ))} \\
 &=-\frac{1}{2}\int_{{\mathbb R}^3} {\theta _{r} u_i^2 (\partial _{x_1 }
u_1 +\partial _{x_2 } u_2 +\partial _{x_3 } u_3 )}
-\frac{1}{2}\int_{{\mathbb R}^3} {u_i^2 (u_1 \partial _{x_1 } \theta
_{r} +u_2 \partial _{x_2
} \theta _{r} +u_3 \partial _{x_3 } \theta _{r} )} \\
 &=-\frac{1}{2}\int_{{\mathbb R}^3} {u_i^2 (u_1 \partial _{x_1 } \theta _{r} +u_2 \partial _{x_2 } \theta _{r} +u_3 \partial _{x_3 }
\theta _{r} )} ,\quad \quad i=1,2,3
\end{split}
\end{equation*}
Taking $ r\to +\infty $ we get
\[
\int_{{\mathbb R}^3} {u_i (u_1 \partial _{x_1 } u_i +u_2 \partial
_{x_2 } u_i +u_3 \partial _{x_3 } u_i )} =0,\quad \quad i=1,2,3
\]
in the same way,
\[
\int_{{\mathbb R}^3} {(u_1 \partial _{x_1 } p+u_2 \partial _{x_2 }
p+u_3
\partial _{x_3 } p)} =0
\]
and
\[
\int_{{\mathbb R}^3} {u_i \Delta u_i } =\int_{{\mathbb R}^3} {u_i
(\partial _{x_1 }^2 u_i +\partial _{x_2 }^2 u_i +\partial _{x_3 }^2
u_i )} =-\int_{{\mathbb R}^3} {((\partial _{x_1 } u_i )^2+(\partial
_{x_2 } u_i )^2+(\partial _{x_3 } u_i )^2)}
\]
then,
\begin{equation*}
\begin{split}
 &\int_{{\mathbb R}^3} {u_1 \partial _t \,u_1 } +\int_{{\mathbb R}^3} {u_1 (u_1
\partial _{x_1 } u_1 +u_2 \partial _{x_2 } u_1 +u_3 \partial _{x_3 } u_1 )}
+\int_{{\mathbb R}^3} {u_1 \partial _{x_1 } p} =\int_{{\mathbb R}^3}
{u_1 \Delta
u_1 } \\
 &\int_{{\mathbb R}^3} {u_2 \partial _t \,u_2 } +\int_{{\mathbb R}^3} {u_2 (u_1
\partial _{x_1 } u_2 +u_2 \partial _{x_2 } u_2 +u_3 \partial _{x_3 } u_2 )}
+\int_{{\mathbb R}^3} {u_2 \partial _{x_2 } p} =\int_{{\mathbb R}^3}
{u_2 \Delta u_2
} \\
 &\int_{{\mathbb R}^3} {u_3 \partial _t \,u_3 } +\int_{{\mathbb R}^3} {u_3 (u_1
\partial _{x_1 } u_3 +u_2 \partial _{x_2 } u_3 +u_3 \partial _{x_3 } u_3 )}
+\int_{{\mathbb R}^3} {u_3 \partial _{x_3 } p} =\int_{{\mathbb R}^3}
{u_3 \Delta
u_3 } \\
\end{split}
\end{equation*}
so that
\begin{equation*}
\begin{split}
 &\frac{1}{2}\partial _t \;\int_{{\mathbb R}^3} {(u_1^2 +u_2^2 +u_3^2 )}
+\;\int_{{\mathbb R}^3} {((\partial _{x_1 } u_1 )^2+(\partial _{x_2
} u_1
)^2+(\partial _{x_3 } u_1 )^2+} \\
 &+(\partial _{x_1 } u_2 )^2+(\partial _{x_2 } u_2 )^2+(\partial _{x_3 } u_2
)^2+(\partial _{x_1 } u_3 )^2+(\partial _{x_2 } u_3 )^2+(\partial _{x_3 }
u_3 )^2)=0 \\
\end{split}
\end{equation*}
it follows that
\begin{equation*}
\begin{split}
 &\int_{{\mathbb R}^3} {(u_1^2 +u_2^2 +u_3^2 )} +2\;\int_0^T {(\,\left\| {\nabla
u_1 } \right\|_{L^2({\mathbb R}^3)}^2 +} \left\| {\nabla u_2 }
\right\|_{L^2({\mathbb R}^3)}^2 +\left\| {\nabla u_3 }
\right\|_{L^2({\mathbb
R}^3)}^2 ) \\
 &\quad \quad =\int_{{\mathbb R}^3} {(u_{10}^2 +u_{20}^2 +u_{30}^2 )} \\
\end{split}
\end{equation*}
Hence from (3) we have
\begin{equation}
\begin{split}
 &\mathop {\sup }\limits_{t\in (0,T)} \;\;\int_{{\mathbb R}^3} {(u_1^2 +u_2^2
+u_3^2 )} <+\infty \\
 &\int_0^T {(\,\left\| {\nabla u_1 } \right\|_{L^2({\mathbb R}^3)}^2 +} \left\|
{\nabla u_2 } \right\|_{L^2({\mathbb R}^3)}^2 +\left\| {\nabla u_3 }
\right\|_{L^2({\mathbb R}^3)}^2 )<+\infty \\
\end{split}
\end{equation}
Above $u$ can be interpreted as the Galerkin approximation of the
solution, but (6) are also true for the solution of problem (1).
\\

The rest of sections are arranged as follows : In section 2 and 3, we introduce a set of auxiliary problems and prove the uniform boundedness and the existence of their solutions in $L^{\infty}(0,T;L^2({\mathbb R}^3))$. Then it is shown that the solutions of the auxiliary problems converge to that of Naiver-Stokes equation with vorticity-velocity form, which also belongs to $L^{\infty}(0,T;L^2({\mathbb R}^3))$. Final section will present the solution of Navier-Stokes equation with velocity-pressure form belongs to $L^{\infty}(0,T; H^2({\mathbb R}^3))$.
\\
\\
\\

\textbf{2. Auxiliary Problems}

For the 3D regularity, we only need to prove that the vorticity in (2)
belongs to $L^\infty (0,T;L^2({\mathbb R}^3))$.
\\

Given a partition with respect to $t$ as follows:
\[
0=t_0 <t_1 <t_2 <\cdots <t_{k-1} <t_k <\cdots <t_N =T
\]
On each $t\in (t_{k-1} ,\;t_k )$, we introduce an auxiliary problem:

\begin{equation}
\begin{split}
 &\partial _t \tilde {\omega }_1 \,+ \overline {u}_1^k \partial _{x_1 } \overline{\overline
{\omega }}_1^k + \overline {u}_2^k \partial _{x_2 }
\overline{\overline {\omega }}_1^k + \overline {u}_3^k \partial
_{x_3 } \overline{\overline {\omega }}_1^k - \overline{\overline
{\omega }}_1^k
\partial _{x_1 } \overline {u}_1^k - \overline{\overline {\omega }}_2^k
\partial _{x_2 } \overline {u}_1^k - \overline{\overline {\omega }}_3^k \partial _{x_3 }
\overline {u}_1^k + \partial _{x_1 } q = \Delta
\tilde {\omega }_1 \\
 &\partial _t \tilde {\omega }_2 + \overline {u}_1^k \partial _{x_1 } \overline{\overline {\omega
}}_2^k + \overline {u}_2^k \partial _{x_2 } \overline{\overline
{\omega }}_2^k + \overline {u}_3^k
\partial _{x_3 } \overline{\overline {\omega }}_2^k - \overline{\overline {\omega }}_1^k \partial _{x_1 }
\overline {u}_2^k - \overline{\overline {\omega }}_2^k \partial
_{x_2 } \overline {u}_2^k - \overline{\overline {\omega }}_3^k
\partial _{x_3 } \overline {u}_2^k + \partial _{x_2 } q = \Delta
\tilde
{\omega }_2 \\
 &\partial _t \tilde {\omega }_3 + \overline {u}_1^k \partial _{x_1 } \overline{\overline {\omega
}}_3^k + \overline {u}_2^k \partial _{x_2 } \overline{\overline
{\omega }}_3^k + \overline {u}_3^k
\partial _{x_3 } \overline{\overline {\omega }}_3^k - \overline{\overline {\omega }}_1^k \partial _{x_1 }
\overline {u}_3^k - \overline{\overline {\omega }}_2^k \partial
_{x_2 } \overline {u}_3^k - \overline{\overline {\omega }}_3^k
\partial _{x_3 } \overline {u}_3^k +
\partial _{x_3 } q = \Delta \tilde
{\omega }_3 \\
\end{split}
\end{equation}
where the initial value is assumed to be $\tilde {\omega }_i (x,t_{k-1}
)=\tilde {\omega }_i^{k-1}(x) $, \;$\tilde {\omega }_i
(x,0 ) = {\omega }_{i0} (x), \,i=1,2,3 $, \; and
\[
 \overline {u}_i^k (x) =
\frac{1}{\Delta t_k }\int_{t_{k-1} }^{t_k } {u_i (x,t)dt}
\]
and
\[
\overline {\omega }_i^k (x) = \frac{1}{\Delta t_k }\int_{t_{k-1}
}^{t_k } {\tilde {\omega }_i (x,t)dt},\quad \quad i=1,2,3
\]

In addition, let $\varepsilon > 0$, we construct a mollifier
$J_\varepsilon \in C_0^\infty (\mathbb R^3)$ such that

\qquad i) $J_\varepsilon (x) \ge 0,\;\;x \in \mathbb R^3$,

\qquad ii) $J_\varepsilon (x) = 0$ if $\;\left| x \right| \ge
\varepsilon $, and

\qquad iii) $\int_{\mathbb R^3} {J_\varepsilon (x)\, dx} = 1$.
\\

\noindent then a convolution is defined as
\[
\overline{\overline {\omega }} _i^k (x) = J_\varepsilon * \overline
{\omega } _i^k (x) = \int_{\mathbb R^3} {J_\varepsilon (x - y)\,
\overline {\omega } _i^k (y)\,dy}
\]

Similarly we can set the incompressible condition :
\begin{equation*}
\begin{split}
 &\partial _{x_1 } u_1 +\partial _{x_2 } u_2 +\partial _{x_3 } u_3 =0 \\
 &\partial _{x_1 } \tilde \omega _1 + \partial _{x_2 } \tilde \omega _2 + \partial _{x_3 }
\tilde \omega _3 = 0 \\
\end{split}
\end{equation*}

It is easy to check that
\begin{equation*}
\begin{split}
 &\partial _{x_1 } u_1 +\partial _{x_2 } u_2 +\partial _{x_3 } u_3 = 0 \quad
\,\Rightarrow \quad \partial _{x_1 } \overline {u}_1^k +\partial
_{x_2 } \overline
{u}_2^k +\partial _{x_3 } \overline {u}_3^k =0 \\
 &\partial _{x_1 } \tilde {\omega }_1 + \partial _{x_2 } \tilde {\omega }_2
+\partial _{x_3 } \tilde {\omega }_3 = 0 \quad \Rightarrow \quad
\partial _{x_1 } \overline{\overline {\omega }}_1^k +\partial _{x_2 } \overline{\overline
{\omega }}_2^k +\partial
_{x_3 } \overline{\overline {\omega }}_3^k =0 \\
\end{split}
\end{equation*}

In the section 3, by means of the Galerkin method and the
compactness imbedding theorem, we can prove the local existences of
the weak solutions of these systems for each $(t_{k-1} ,\;t_k )$
being small enough. Below we also interpret $\tilde \omega$ as the
Galerkin approximation of the solution of the problems (7), and
first prove that  $\tilde \omega, t \in (0,T)$ belong to
$L^{\infty}(0,T;L^2({\mathbb R}^3))$. In section 4, an approach of
approximation is used to assert that the solution of (2) also
belongs to $L^\infty (0,T;L^2({\mathbb R}^3))$.
\\

Since
\begin{equation*}
\begin{split}
 &\int_{{\mathbb R}^3} {\theta _{r} \;\; [\,\tilde {\omega }_1 (\overline {u}_1^k
\partial _{x_1 } \overline{\overline {\omega }}_1^k + \overline {u}_2^k \partial _{x_2 } \overline{\overline {\omega
}}_1^k + \overline {u}_3^k \partial _{x_3 } \overline{\overline {\omega }}_1^k )} \\
 &\;\;\qquad + \tilde {\omega }_2 (\overline {u}_1^k \partial _{x_1 } \overline{\overline
{\omega }}_2^k + \overline {u}_2^k \partial _{x_2 }
\overline{\overline {\omega }}_2^k + \overline
{u}_3^k \partial _{x_3 } \overline{\overline {\omega }}_2^k ) \\
 &\;\;\qquad + \tilde {\omega }_3 (\overline {u}_1^k \partial _{x_1 } \overline{\overline
{\omega }}_3^k + \overline {u}_2^k \partial _{x_2 }
\overline{\overline {\omega }}_3^k + \overline
{u}_3^k \partial _{x_3 } \overline{\overline {\omega }}_3^k )\,] \\
 &= -\int_{{\mathbb R}^3} {[\,\overline{\overline {\omega }}_1^k \partial _{x_1 } (\theta _{r} \tilde {\omega }_1 \overline {u}_1^k ) + \overline{\overline {\omega }}_1^k \partial
_{x_2 } (\theta _{r} \tilde {\omega }_1 \overline {u}_2^k )+
\overline{\overline {\omega }}_1^k
\partial _{x_3 } (\theta _{r} \tilde {\omega }_1 \overline {u}_3^k )} \\
 &\quad \qquad + \overline{\overline {\omega }}_2^k \partial _{x_1 } (\theta _{r}
\tilde {\omega }_2 \overline {u}_1^k ) + \overline{\overline {\omega
}}_2^k
\partial _{x_2 } (\theta _{r} \tilde {\omega }_2 \overline {u}_2^k ) + \overline{\overline
{\omega }}_2^k
\partial _{x_3 } (\theta _{r} \tilde {\omega }_2 \overline {u}_3^k ) \\
 &\quad \qquad + \overline{\overline {\omega }}_3^k \partial _{x_1 } (\theta _{r}
\tilde {\omega }_3 \overline {u}_1^k ) + \overline{\overline {\omega
}}_3^k
\partial _{x_2 } (\theta _{r} \tilde {\omega }_3 \overline {u}_2^k )
+ \overline{\overline {\omega }}_3^k
\partial _{x_3 } (\theta _{r} \tilde {\omega }_3 \overline {u}_3^k )\,] \\
 &= -\int_{{\mathbb R}^3} {(\overline{\overline {\omega }}_1^k \tilde {\omega }_1 \overline {u}_1^k
\partial _{x_1 } \theta _{r} + \overline{\overline {\omega }}_1^k \theta _{r} \overline {u}_1^k \partial _{x_1 } \tilde {\omega }_1 + \overline{\overline {\omega }}_1^k
\theta _{r} \tilde {\omega }_1 \partial _{x_1 } \overline {u}_1^k } \\
 &\quad \qquad + \overline{\overline {\omega }}_1^k \tilde {\omega }_1 \overline {u}_2^k \partial
_{x_2 } \theta _{r} + \overline{\overline {\omega }}_1^k \theta _{r}
\overline {u}_2^k
\partial _{x_2 } \tilde {\omega }_1 + \overline{\overline {\omega }}_1^k \theta
_{r} \tilde {\omega }_1 \partial _{x_2 } \overline {u}_2^k \\
 &\quad \qquad + \overline{\overline {\omega }}_1^k \tilde {\omega }_1 \overline {u}_3^k \partial
_{x_3 } \theta _{r} + \overline{\overline {\omega }}_1^k \theta _{r}
\overline {u}_3^k
\partial _{x_3 } \tilde {\omega }_1 + \overline{\overline {\omega }}_1^k \theta
_{r} \tilde {\omega }_1 \partial _{x_3 } \overline {u}_3^k \\
 &\quad \qquad + \overline {\omega }_2^k \tilde {\omega }_2 \overline {u}_1^k \partial
_{x_1 } \theta _{r} + \overline{\overline {\omega }}_2^k \theta _{r}
\overline {u}_1^k
\partial _{x_1 } \tilde {\omega }_2 + \overline{\overline {\omega }}_2^k \theta
_{r} \tilde {\omega }_2 \partial _{x_1 } \overline {u}_1^k \\
 &\quad \qquad + \overline{\overline {\omega }}_2^k \tilde {\omega }_2 \overline {u}_2^k \partial
_{x_2 } \theta _{r} + \overline{\overline {\omega }}_2^k \theta _{r}
\overline {u}_2^k
\partial _{x_2 } \tilde {\omega }_2 + \overline{\overline {\omega }}_2^k \theta
_{r} \tilde {\omega }_2 \partial _{x_2 } \overline {u}_2^k \\
 &\quad \qquad + \overline{\overline {\omega }}_2^k \tilde {\omega }_2 \overline {u}_3^k \partial
_{x_3 } \theta _{r} + \overline{\overline {\omega }}_2^k \theta _{r}
\overline {u}_3^k
\partial _{x_3 } \tilde {\omega }_2 + \overline{\overline {\omega }}_2^k \theta
_{r} \tilde {\omega }_2 \partial _{x_3 } \overline {u}_3^k \\
 &\quad \qquad + \overline{\overline {\omega }}_3^k \tilde {\omega }_3 \overline {u}_1^k \partial
_{x_1 } \theta _{r} + \overline{\overline {\omega }}_3^k \theta _{r}
\overline {u}_1^k
\partial _{x_1 } \tilde {\omega }_3 + \overline{\overline {\omega }}_3^k \theta
_{r} \tilde {\omega }_3 \partial _{x_1 } \overline {u}_1^k \\
 &\quad \qquad + \overline{\overline {\omega }}_3^k \tilde {\omega }_3 \overline {u}_2^k \partial
_{x_2 } \theta _{r} + \overline{\overline {\omega }}_3^k \theta _{r}
\overline {u}_2^k
\partial _{x_2 } \tilde {\omega }_3 + \overline{\overline {\omega }}_3^k \theta
_{r} \tilde {\omega }_3 \partial _{x_2 } \overline {u}_2^k \\
 &\quad \qquad + \overline{\overline {\omega }}_3^k \tilde {\omega }_3 \overline {u}_3^k \partial
_{x_3 } \theta _{r} + \overline{\overline {\omega }}_3^k \theta _{r}
\overline {u}_3^k
\partial _{x_3 } \tilde {\omega }_3 + \overline{\overline {\omega }}_3^k \theta
_{r} \tilde {\omega }_3 \partial _{x_3 } \overline {u}_3^k ) \\
 &\\
 &= -\int_{{\mathbb R}^3} {[\,\theta _{r} (\overline{\overline {\omega }}_1^k \overline {u}_1^k
\partial _{x_1 } \tilde {\omega }_1 + \overline{\overline {\omega }}_1^k \overline {u}_2^k
\partial _{x_2 } \tilde {\omega }_1 + \overline{\overline {\omega }}_1^k \overline {u}_3^k
\partial _{x_3 } \tilde {\omega }_1 } \\
 &\qquad \quad \quad \; + \overline{\overline {\omega }}_2^k \overline {u}_1^k \partial _{x_1 }
\tilde {\omega }_2 + \overline{\overline {\omega }}_2^k \overline
{u}_2^k
\partial _{x_2 } \tilde {\omega }_2 + \overline{\overline {\omega }}_2^k \overline
{u}_3^k
\partial _{x_3 } \tilde {\omega
}_2 \\
 &\qquad \quad \quad \; + \overline{\overline {\omega }}_3^k \overline {u}_1^k \partial _{x_1 }
\tilde {\omega }_3 + \overline{\overline {\omega }}_3^k \overline
{u}_2^k
\partial _{x_2 } \tilde {\omega }_3 + \overline{\overline {\omega }}_3^k \overline
{u}_3^k
\partial _{x_3 } \tilde
{\omega }_3 ) \\
 &\qquad \quad \;   +\; (\overline{\overline {\omega }}_1^k \tilde {\omega }_1 \overline {u}_1^k \partial _{x_1 }
\theta _{r} + \overline{\overline {\omega }}_1^k \tilde {\omega }_1
\overline {u}_2^k
\partial _{x_2 } \theta _{r} + \overline{\overline {\omega }}_1^k \tilde {\omega }_1
\overline {u}_3^k \partial _{x_3 } \theta _{r} \\
 &\qquad \quad \;   +\; \overline{\overline {\omega }}_2^k \tilde {\omega }_2 \overline {u}_1^k \partial
_{x_1 } \theta _{r} + \overline{\overline {\omega }}_2^k \tilde
{\omega }_2 \overline {u}_2^k \partial _{x_2 } \theta _{r} +
\overline{\overline {\omega }}_2^k \tilde
{\omega }_2 \overline {u}_3^k \partial _{x_3 } \theta _{r} \\
 &\qquad \quad \;   +\; \overline{\overline {\omega }}_3^k \tilde {\omega }_3 \overline {u}_1^k \partial
_{x_1 } \theta _{r} + \overline{\overline {\omega }}_3^k \tilde
{\omega }_3 \overline {u}_2^k \partial _{x_2 } \theta _{r} +
\overline{\overline {\omega }}_3^k \tilde
{\omega }_3 \overline {u}_3^k \partial _{x_3 } \theta _{r} )\,] \\
\end{split}
\end{equation*}
\\
Let $ r\to +\infty $ we get
\begin{equation*}
\begin{split}
 &\int_{{\mathbb R}^3} {[\,\tilde {\omega }_1 (\overline {u}_1^k \partial _{x_1 } \overline{\overline
{\omega }}_1^k +\overline {u}_2^k \partial _{x_2 }
\overline{\overline {\omega}} _1^k + \overline {u}_3^k
\partial _{x_3 } \overline{\overline {\omega}} _1^k )} \\
 &\;\;\;\,+ \tilde {\omega }_2 (\overline {u}_1^k \partial _{x_1 } \overline{\overline {\omega
}}_2^k + \overline {u}_2^k \partial _{x_2 } \overline{\overline
{\omega }}_2^k + \overline {u}_3^k
\partial _{x_3 } \overline{\overline {\omega }}_2^k ) \\
 &\;\;\;\,+ \tilde {\omega }_3 (\overline {u}_1^k \partial _{x_1 } \overline{\overline {\omega
}}_3^k + \overline {u}_2^k \partial _{x_2 } \overline{\overline
{\omega }}_3^k + \overline {u}_3^k
\partial _{x_3 } \overline{\overline {\omega }}_3^k )\,] \\
 &= -\int_{{\mathbb R}^3} {(\overline{\overline {\omega }}_1^k \overline {u}_1^k \partial _{x_1 }
\tilde {\omega }_1 + \overline{\overline {\omega }}_1^k \overline
{u}_2^k
\partial _{x_2 } \tilde {\omega }_1 + \overline{\overline {\omega }}_1^k \overline
{u}_3^k
\partial _{x_3 } \tilde
{\omega }_1 } \;\, \\
 &\quad \qquad  + \overline{\overline {\omega }}_2^k \overline {u}_1^k \partial _{x_1 } \tilde
{\omega }_2 + \overline{\overline {\omega }}_2^k \overline {u}_2^k
\partial _{x_2 } \tilde {\omega
}_2 + \overline{\overline {\omega }}_2^k \overline {u}_3^k \partial _{x_3 } \tilde {\omega }_2 \\
 &\quad \qquad  + \overline{\overline {\omega }}_3^k \overline {u}_1^k \partial _{x_1 } \tilde
{\omega }_3 + \overline{\overline {\omega }}_3^k \overline {u}_2^k
\partial _{x_2 } \tilde {\omega }_3 + \overline{\overline {\omega }}_3^k
\overline {u}_3^k
\partial _{x_3 } \tilde {\omega }_3
\,) \\
\end{split}
\end{equation*}
Similarly,
\begin{equation*}
\begin{split}
 &\int_{{\mathbb R}^3} {[\,\tilde {\omega }_1 (\overline{\overline {\omega }}_1^k \partial _{x_1 }
\overline {u}_1^k + \overline{\overline {\omega }}_2^k \partial
_{x_2 } \overline {u}_1^k + \overline{\overline {\omega
}}_3^k \partial _{x_3 } \overline {u}_1^k } ) \\
 &\;\,\;\,+ \tilde {\omega }_2 (\overline{\overline {\omega }}_1^k \partial _{x_1 } \overline
{u}_2^k + \overline{\overline {\omega }}_2^k \partial _{x_2 }
\overline {u}_2^k + \overline{\overline {\omega
}}_3^k \partial _{x_3 } \overline {u}_2^k ) \\
 &\;\;\,\, + \tilde {\omega }_3 (\overline{\overline {\omega }}_1^k \partial _{x_1 } \overline
{u}_3^k + \overline{\overline {\omega }}_2^k \partial _{x_2 }
\overline {u}_3^k + \overline{\overline {\omega
}}_3^k \partial _{x_3 } \overline {u}_3^k )\,] \\
 &= -\int_{{\mathbb R}^3} {(\overline{\overline {\omega }}_1^k \overline {u}_1^k \partial _{x_1 }
\tilde {\omega }_1 + \overline{\overline {\omega }}_2^k \overline
{u}_1^k
\partial _{x_2 } \tilde {\omega }_1 + \overline{\overline {\omega }}_3^k \overline
{u}_1^k
\partial _{x_3 }
\tilde {\omega }_1 } \\
 &\quad \qquad + \overline{\overline {\omega }}_1^k \overline {u}_2^k \partial _{x_1 } \tilde
{\omega }_2 + \overline{\overline {\omega }}_2^k \overline {u}_2^k
\partial _{x_2 } \tilde {\omega }_2 + \overline{\overline {\omega }}_3^k
\overline {u}_2^k
\partial _{x_3 } \tilde
{\omega }_2 \\
 &\quad \qquad + \overline{\overline {\omega }}_1^k \overline {u}_3^k \partial _{x_1 } \tilde
{\omega }_3 + \overline{\overline {\omega }}_2^k \overline {u}_3^k
\partial _{x_2 } \tilde {\omega }_3 + \overline{\overline {\omega }}_3^k
\overline {u}_3^k
\partial _{x_3 } \tilde
{\omega }_3 ) \\
\end{split}
\end{equation*}
and
\[
\int_{{\mathbb R}^3} {(\tilde {\omega }_1 \partial _{x_1 } q +
\tilde {\omega }_2
\partial _{x_2 } q + \tilde {\omega }_3 \partial _{x_3 } q)} =0
\]
furthermore,
\[
\int_{{\mathbb R}^3} {\tilde {\omega }_i \Delta \tilde {\omega }_i }
=\int_{{\mathbb R}^3} {\tilde {\omega }_i (\partial _{x_1 }^2 \tilde
{\omega }_i +\partial _{x_2 }^2 \tilde {\omega }_i +\partial _{x_3
}^2 \tilde {\omega }_i )} =-\int_{{\mathbb R}^3} {((\partial _{x_1 }
\tilde {\omega }_i )^2+(\partial _{x_2 } \tilde {\omega }_i
)^2+(\partial _{x_3 } \tilde {\omega }_i )^2)}
\]
\\

Thus from (7) we have
\begin{equation*}
\begin{split}
 &\int_{{\mathbb R}^3} {\tilde {\omega }_1 \partial _t \tilde {\omega }_1 }
\;\,+ \int_{{\mathbb R}^3} {\tilde {\omega }_1 (\overline {u}_1^k
\partial _{x_1 } \overline{\overline {\omega }}_1^k + \overline {u}_2^k \partial _{x_2 }
\overline{\overline {\omega }}_1^k + \overline
{u}_3^k \partial _{x_3 } \overline{\overline {\omega }}_1^k )} \\
 &\quad \quad \quad \quad \quad \quad \quad \quad -\int_{{\mathbb R}^3} {\tilde
{\omega }_1 (\overline{\overline {\omega }}_1^k \partial _{x_1 }
\overline {u}_1^k + \overline{\overline {\omega }}_2^k \partial
_{x_2 } \overline {u}_1^k + \overline{\overline {\omega }}_3^k
\partial _{x_3 } \overline {u}_1^k )} \,+ \int_{{\mathbb R}^3} {\tilde {\omega
}_1 \partial _{x_1 } q}
= \int_{{\mathbb R}^3} {\tilde {\omega }_1 \Delta \tilde {\omega }_1 } \\
 &\int_{{\mathbb R}^3} {\tilde {\omega }_2 \partial _t \tilde {\omega }_2 }
+\int_{{\mathbb R}^3} {\tilde {\omega }_2 (\overline {u}_1^k
\partial _{x_1 } \overline{\overline {\omega }}_2^k + \overline {u}_2^k \partial _{x_2 }
\overline{\overline {\omega }}_2^k + \overline
{u}_3^k \partial _{x_3 } \overline{\overline {\omega }}_2^k )} \\
 &\quad \quad \quad \quad \quad \quad \quad \quad -\int_{{\mathbb R}^3} {\tilde
{\omega }_2 (\overline{\overline {\omega }}_1^k \partial _{x_1 }
\overline {u}_2^k + \overline{\overline {\omega }}_2^k \partial
_{x_2 } \overline {u}_2^k + \overline{\overline {\omega }}_3^k
\partial _{x_3 } \overline {u}_2^k )} + \int_{{\mathbb R}^3} {\tilde {\omega
}_2 \partial _{x_2 } q}
=\int_{{\mathbb R}^3} {\tilde {\omega }_2 \Delta \tilde {\omega }_2 } \\
\end{split}
\end{equation*}
\begin{equation*}
\begin{split}
 &\int_{{\mathbb R}^3} {\tilde {\omega }_3 \partial _t \tilde {\omega }_3 }
+ \int_{{\mathbb R}^3} {\tilde {\omega }_3 (\overline {u}_1^k
\partial _{x_1 } \overline{\overline {\omega }}_3^k + \overline {u}_2^k \partial _{x_2 }
\overline{\overline {\omega }}_3^k + \overline
{u}_3^k \partial _{x_3 } \overline{\overline {\omega }}_3^k )} \\
 &\quad \quad \quad \quad \quad \quad \quad \quad -\int_{{\mathbb R}^3} {\tilde
{\omega }_3 (\overline{\overline {\omega }}_1^k \partial _{x_1 }
\overline {u}_3^k + \overline{\overline {\omega }}_2^k \partial
_{x_2 } \overline {u}_3^k + \overline{\overline {\omega }}_3^k
\partial _{x_3 } \overline {u}_3^k )} + \int_{{\mathbb R}^3} {\tilde {\omega
}_3 \partial _{x_3 } q}
=\int_{{\mathbb R}^3} {\tilde {\omega }_3 \Delta \tilde {\omega }_3 } \\
\end{split}
\end{equation*}
so that
\begin{equation*}
\begin{split}
 &\frac{1}{2}\partial _t \int_{{\mathbb R}^3} {(\tilde {\omega }_1^2 +\tilde
{\omega }_2^2 +\tilde {\omega }_3^2 )} +\,\,\int_{{\mathbb R}^3}
{[\,(\partial _{x_1 } \tilde {\omega }_1 )^2 + (\partial _{x_2 }
\tilde {\omega }_1
)^2+(\partial _{x_3 } \tilde {\omega }_1 )^2 } \\
 &\qquad \quad \quad \quad \quad \quad \quad \quad \quad \quad
\;\,\;\;\;\, + (\partial _{x_1 } \tilde {\omega }_2 )^2 + (\partial
_{x_2 } \tilde
{\omega }_2 )^2 + (\partial _{x_3 } \tilde {\omega }_2 )^2 \\
 &\qquad \quad \quad \quad \quad \quad \quad \quad \quad \quad
\;\,\;\;\;\, + (\partial _{x_1 } \tilde {\omega }_3 )^2 + (\partial
_{x_2 } \tilde
{\omega }_3 )^2 + (\partial _{x_3 } \tilde {\omega }_3 )^2] \\
 &\quad \quad \quad -\int_{{\mathbb R}^3} {(\overline{\overline {\omega }}_1^k \overline {u}_1^k
\partial _{x_1 } \tilde {\omega }_1 + \overline{\overline {\omega }}_1^k \overline {u}_2^k
\partial _{x_2 } \tilde {\omega }_1 + \overline{\overline {\omega }}_1^k \overline {u}_3^k
\partial _{x_3 } \tilde {\omega }_1 } \\
 &\quad \quad \quad \,\quad \;\;\; + \overline{\overline {\omega }}_2^k \overline {u}_1^k \partial
_{x_1 } \tilde {\omega }_2 + \overline{\overline {\omega }}_2^k
\overline {u}_2^k
\partial _{x_2 } \tilde {\omega }_2 + \overline{\overline {\omega }}_2^k \overline
{u}_3^k
\partial _{x_3 } \tilde
{\omega }_2 \\
 &\quad \quad \quad \,\quad \;\;\; + \overline{\overline {\omega }}_3^k \overline {u}_1^k \partial
_{x_1 } \tilde {\omega }_3 + \overline{\overline {\omega }}_3^k
\overline {u}_2^k
\partial _{x_2 } \tilde {\omega }_3 + \overline{\overline {\omega }}_3^k \overline
{u}_3^k \partial _{x_3 }
\tilde {\omega }_3 ) \\
 &\quad \quad \quad +\int_{{\mathbb R}^3} {(\overline{\overline {\omega }}_1^k \overline {u}_1^k
\partial _{x_1 } \tilde {\omega }_1 + \overline{\overline {\omega }}_2^k \overline {u}_1^k
\partial _{x_2 } \tilde {\omega }_1 + \overline{\overline {\omega }}_3^k \overline {u}_1^k
\partial _{x_3 } \tilde {\omega }_1 } \\
 &\quad \quad \quad \,\quad \;\;\; + \overline{\overline {\omega }}_1^k \overline {u}_2^k \partial
_{x_1 } \tilde {\omega }_2 + \overline{\overline {\omega }}_2^k
\overline {u}_2^k
\partial _{x_2 } \tilde {\omega }_2 + \overline{\overline {\omega }}_3^k \overline {u}_2^k
\partial _{x_3 } \tilde
{\omega }_2 \\
 &\quad \quad \quad \,\quad \;\;\; + \overline{\overline {\omega }}_1^k \overline {u}_3^k \partial
_{x_1 } \tilde {\omega }_3 + \overline{\overline {\omega }}_2^k
\overline {u}_3^k
\partial _{x_2 } \tilde {\omega }_3 + \overline{\overline {\omega }}_3^k \overline
{u}_3^k \partial _{x_3 }
\tilde {\omega }_3 )=0 \\
\end{split}
\end{equation*}
\\

By using Young inequality: $uv\le \frac{1}{4}u^2+v^2$, it follows
that
\begin{equation*}
\begin{split}
 &\partial _t \int_{{\mathbb R}^3} {(\tilde {\omega }_1^2 + \tilde {\omega }_2^2 + \tilde
{\omega }_3^2 )} \;\,+\,\;2\; {\int_{{\mathbb R}^3} {\;[\,(\partial
_{x_1 } \tilde {\omega }_1 )^2\; + (\partial _{x_2 } \tilde
{\omega }_1 )^2 + (\partial _{x_3 } \tilde {\omega }_1 )^2} } \\
 &\qquad \qquad \qquad \quad \quad \quad \quad \quad \quad
\;\;\;\;\; +(\partial _{x_1 } \tilde {\omega }_2 )^2+(\partial _{x_2
}
\tilde {\omega }_2 )^2+(\partial _{x_3 } \tilde {\omega }_2 )^2 \\
 &\qquad \qquad \qquad \quad \quad \quad \quad \quad \quad
\;\;\;\;\; +(\partial _{x_1 } \tilde {\omega }_3 )^2+(\partial _{x_2
}
\tilde {\omega }_3 )^2 + (\partial _{x_3 } \tilde {\omega }_3 )^2] \\
 &\quad \quad \quad \quad \quad \quad \quad  \le 2 \; {\int_{{\mathbb R}^3} {(\overline{\overline {\omega }}_1^{k^2} \overline {u}_1^{k^2}
+ \overline{\overline {\omega
}}_1^{k^2} \overline {u}_2^{k^2} + \overline{\overline {\omega }}_1^{k^2} \overline {u}_3^{k^2} } } \\
 &\qquad \qquad \qquad \quad \;
\quad \quad \;\; + \overline{\overline {\omega }}_2^{k^2} \overline
{u}_1^{k^2} + \overline{\overline {\omega
}}_2^{k^2} \overline {u}_2^{k^2} + \overline{\overline {\omega }}_2^{k^2} \overline {u}_3^{k^2} \\
 &\qquad \qquad \qquad \quad \;
\quad \quad \;\; + \overline{\overline {\omega }}_3^{k^2} \overline
{u}_1^{k^2} + \overline{\overline {\omega
}}_3^{k^2} \overline {u}_2^{k^2} + \overline{\overline {\omega }}_3^{k^2} \overline {u}_3^{k^2} ) \\
 &\quad \quad \quad \quad \quad \quad \quad \quad \quad \quad
\;\;\;\,+\;\,2 \; {\int_{{\mathbb R}^3} {(\overline{\overline
{\omega }}_1^{k^2} \overline {u}_1^{k^2} + \overline{\overline
{\omega }}_2^{k^2} \overline {u}_1^{k^2} + \overline{\overline
{\omega}
}_3^{k^2} \overline {u}_1^{k^2} } } \\
 &\qquad \qquad \quad \quad \quad \quad \quad \quad \quad \quad \,
\quad \quad + \overline{\overline {\omega }}_1^{k^2} \overline
{u}_2^{k^2} + \overline{\overline {\omega }}_2^{k^2}
\overline {u}_2^{k^2} + \overline{\overline {\omega }}_3^{k^2} \overline {u}_2^{k^2} \\
 &\qquad \qquad \quad \quad \quad \quad \quad \quad \quad \quad  \,
\quad \quad + \overline{\overline {\omega }}_1^{k^2} \overline
{u}_3^{k^2} + \overline{\overline {\omega
}}_2^{k^2} \overline {u}_3^{k^2} + \overline{\overline {\omega }}_3^{k^2} \overline {u}_3^{k^2} ) \\
 &\quad \quad \quad \quad \qquad \qquad \qquad +  {\int_{{\mathbb
R}^3} {[\,(\partial _{x_1 } \tilde {\omega }_1 )^2 + (\partial _{x_2
} \tilde
{\omega }_1 )^2 + (\partial _{x_3 } \tilde {\omega }_1 )^2} } \\
 &\qquad \qquad \qquad \qquad \qquad \qquad \,  + (\partial _{x_1
} \tilde {\omega }_2 )^2+(\partial _{x_2 } \tilde {\omega }_2 )^2 +
(\partial
_{x_3 } \tilde {\omega }_2 )^2 \\
 &\qquad \qquad \qquad \qquad \qquad \qquad \,  + (\partial _{x_1
} \tilde {\omega }_3 )^2 + (\partial _{x_2 } \tilde {\omega }_3 )^2
+ (\partial
_{x_3 } \tilde {\omega }_3 )^2] \\
\end{split}
\end{equation*}
Thus,
\begin{equation}
\label{eq8}
\begin{split}
 &\partial _t \int_{{\mathbb R}^3} {(\tilde {\omega }_1^2 + \tilde {\omega }_2^2 +\tilde
{\omega }_3^2 )} +  {\int_{{\mathbb R}^3} {[\,(\partial _{x_1 }
\tilde {\omega }_1 )^2 + (\partial _{x_2 } \tilde {\omega }_1 )^2 +
(\partial
_{x_3 } \tilde {\omega }_1 )^2} } \\
 &\quad \qquad \qquad \quad \quad \quad \quad \quad  \;\;\;\;\; + (\partial
_{x_1 } \tilde {\omega }_2 )^2 + (\partial _{x_2 } \tilde {\omega
}_2
)^2+(\partial _{x_3 } \tilde {\omega }_2 )^2 \\
 &\quad \qquad \qquad \quad \quad \quad \quad \quad  \;\;\;\;\; + (\partial
_{x_1 } \tilde {\omega }_3 )^2+(\partial _{x_2 } \tilde {\omega }_3
)^2 + (\partial _{x_3 } \tilde {\omega }_3 )^2] \\
 &\qquad \le 4\;\{\,\,\left\| {\overline
{u}_1^k } \right\|_{L^2({\mathbb R}^3 )}^2 (\,\left\|
{\overline{\overline {\omega }}_1^k } \right\|_{L^\infty({\mathbb R}^3
)}^2 +\left\| {\overline{\overline {\omega }}_2^k }
\right\|_{L^\infty({\mathbb R}^3 )}^2 +\left\| {\overline{\overline
{\omega }}_3^k }
\right\|_{L^\infty({\mathbb R}^3 )}^2 ) \\
 &\qquad \quad \;\;\;
 + \left\| {\overline {u}_2^k } \right\|_{L^2({\mathbb R}^3
)}^2 (\,\left\| {\overline{\overline {\omega }}_1^k }
\right\|_{L^\infty({\mathbb R}^3 )}^2 + \left\| {\overline{\overline
{\omega }}_2^k } \right\|_{L^\infty({\mathbb R}^3 )}^2 + \left\|
{\overline{\overline {\omega }}_3^k }
\right\|_{L^\infty({\mathbb R}^3 )}^2 ) \\
 &\qquad \quad \;\;\;
 +\left\| {\overline {u}_3^k } \right\|_{L^2({\mathbb R}^3
)}^2 (\,\left\| {\overline{\overline {\omega }}_1^k }
\right\|_{L^\infty({\mathbb R}^3 )}^2 + \left\| {\overline{\overline
{\omega }}_2^k } \right\|_{L^\infty({\mathbb R}^3 )}^2 + \left\|
{\overline{\overline {\omega }}_3^k }
\right\|_{L^\infty({\mathbb R}^3 )}^2 )\,\} \\
 &\qquad = 4\,(\,\left\| {\overline {u}_1^k } \right\|_{L^2({\mathbb R}^3 )}^2
+\left\| {\overline {u}_2^k } \right\|_{L^2({\mathbb R}^3 )}^2 +
\left\| {\overline {u}_3^k } \right\|_{L^2({\mathbb R}^3 )}^2
)\,(\,\left\| {\overline{\overline {\omega }}_1^k }
\right\|_{L^\infty({\mathbb R}^3 )}^2 + \left\| {\overline{\overline
{\omega }}_2^k } \right\|_{L^\infty({\mathbb R}^3 )}^2 + \left\|
{\overline{\overline {\omega }}_3^k }
\right\|_{L^\infty({\mathbb R}^3 )}^2 ) \\
\end{split}
\end{equation}

Note that
\begin{equation*}
\begin{split}
 &\left\| {\overline {u }_i^{k} } \right\|_{L^2({\mathbb R^3} )}^2 =\int_{\mathbb R^3}
{\left( {\frac{1}{\Delta t_{k} }\int_{t_{k-1} }^{t_{k} } {
{u }_i (x,t)dt} } \right)} ^2\le \frac{1}{\Delta t_{k}^2
}\int_{\mathbb R^3} {\Delta
t_{k} \int_{t_{k-1} }^{t_{k} } { {u }_i^2 (x,t)dt} } \\
 &\quad \quad \quad \quad \;\; =\frac{1}{\Delta t_{k} }\int_{t_{k-1} }^{t_{k}
} {\left\| { {u }_i } \right\|_{L^2({\mathbb R^3} )}^2 } \;\; \le
\;\mathop {\sup }\limits_{(t_{k - 1} ,t_k )} \left\| {u_i }
\right\|_{L^2({\mathbb R^3} )}^2  \\
 \end{split}
\end{equation*}
and similarly
\[
 \left\| {\overline \omega_i^{k} } \right\|_{L^2({\mathbb R^3} )}^2
\le \frac{1}{\Delta t_{k} }\int_{t_{k-1} }^{t_{k}
} {\left\| \tilde \omega_i \right\|_{L^2({\mathbb R^3} )}^2 }, \qquad i=1,2,3 \\
\]

In addition, a convolution inequality in [1] is applied to get
\begin{equation*}
\begin{split}
 &\left\| {\overline{\overline \omega} _i^k } \right\|_{L^\infty ({\mathbb R^3} )}^2 =
\left\| {J_\varepsilon * \overline \omega _i^k (x)} \right\|_{L^\infty
({\mathbb R^3}
)}^2 \\
 &\qquad \qquad \le \left\| {J_\varepsilon } \right\|_{L^2({\mathbb R^3} )}^2
\;\left\| {\overline \omega _i^k } \right\|_{L^2({\mathbb R^3} )}^2 \;\; \le \;\frac{1}{\mu_\varepsilon}
\;\mathop {\sup }\limits_{(t_{k - 1} ,t_k )} \left\| {\tilde \omega_i }
\right\|_{L^2({\mathbb R^3} )}^2
\\
\end{split}
\end{equation*}
where the quantity $\mu_\varepsilon \to 0$ as $\varepsilon \to 0$. We need further assuming that $\varepsilon \to 0$ and
\[
\frac{\Delta t_k}{\mu_\varepsilon} \to 0  \qquad \mbox{as}\;\; k \to \infty \;\;\mbox{or}\;\; \Delta t_k \to 0
\]

From (8) we have
\begin{equation*}
\begin{split}
 &\int_{\mathbb R^3} {(\tilde {\omega }_1^2 + \tilde {\omega }_2^2 +\tilde {\omega
}_3^2 )} \;+\; \int_{t_{k-1} }^t {(\,\left\| {\nabla \tilde {\omega
}_1 } \right\|_{L^2({\mathbb R^3} )}^2 +\left\| {\nabla \tilde {\omega }_2
} \right\|_{L^2({\mathbb R^3} )}^2 +\left\| {\nabla \tilde {\omega }_3 }
\right\|_{L^2({\mathbb R^3} )}^2 )} \\
 &\qquad  \le \int_{\mathbb R^3} {(\tilde {\omega }_1^{k-1^2} +\tilde {\omega }_2^{k-1^2} +\tilde
{\omega }_3^{k-1^2} )} \;\; + \\
 &\qquad \qquad  + \;\frac{4\Delta t_{k} }{\mu_\varepsilon} \; \mathop {\sup }\limits_{(t_{k-1} ,t_k )} \left\{
{\left\| {u_1 } \right\|_{L^2({\mathbb R^3} )}^2 +\left\| {u_2 }
\right\|_{L^2({\mathbb R^3} )}^2 +\left\| {u_3 } \right\|_{L^2({\mathbb R^3} )}^2
} \right\}  \cdot \mathop {\sup }\limits_{(t_{k-1} ,t )}\; \int_{\mathbb R^3} {(\tilde {\omega }_1^2 + \tilde {\omega }_2^2 +\tilde {\omega
}_3^2 )}  \\
 \end{split}
\end{equation*}

By (6) we have
\begin{equation*}
\begin{split}
  &\mathop {\sup }\limits_{(t_{k-1} ,t_k )} \left\{
{\left\| {u_1 } \right\|_{L^2({\mathbb R^3} )}^2 +\left\| {u_2 }
\right\|_{L^2({\mathbb R^3} )}^2 +\left\| {u_3 } \right\|_{L^2({\mathbb R^3} )}^2
} \right\}  \\
  &\qquad  \le K_0 = \mathop {\sup }\limits_{t\in (0,T)} \; \int_{\mathbb R^3} {(u_1^2
+u_2^2 +u_3^2 )} + \int_0
^T {\left( {\left\| {\nabla u_1 }
\right\|_{L^2({\mathbb R^3} )}^2 +\left\| {\nabla u_2 }
\right\|_{L^2({\mathbb R^3} )}^2 +\left\| {\nabla u_3 }
\right\|_{L^2({\mathbb R^3} )}^2 } \right)}\; < +\infty   \\
 \end{split}
\end{equation*}
Thus,
\begin{equation*}
\begin{split}
 & {\left( 1-4K_0 \frac{\Delta t_{k} }{\mu_\varepsilon} \right)}\; \mathop {\sup }\limits_{t\in (t_{k-1} ,t_k )} \int_{\mathbb R^3} {(\tilde {\omega
}_1^2 \;+ \tilde {\omega }_2^2 +\tilde {\omega }_3^2 )} \;\; + \\
 &\quad + \int_{t_{k-1}
}^{t_k } {\left( {\left\| {\nabla \tilde {\omega }_1 }
\right\|_{L^2({\mathbb R^3} )}^2 +\left\| {\nabla \tilde {\omega }_2 }
\right\|_{L^2({\mathbb R^3} )}^2 +\left\| {\nabla \tilde {\omega }_3 }
\right\|_{L^2({\mathbb R^3} )}^2 } \right)} \le \int_{\mathbb R^3} {(\tilde {\omega }_1^{k-1^2} +\tilde {\omega }_2^{k-1^2}
+\tilde {\omega }_3^{k-1^2} )}  \\
 \end{split}
\end{equation*}

Now we set
\begin{equation*}
\begin{split}
 & M_0 =\int_{\mathbb R^3} {(\omega _{10}^2 +\omega _{20}^2 +\omega _{30}^2 )}  \\
 & M_k = \mathop {\sup }\limits_{t\in (t_{k-1} ,t_k ) } \;\int_{\mathbb R^3} {(\tilde
{\omega }_1^2 +\tilde {\omega }_2^2 +\tilde {\omega }_3^2 )} \\
 & \delta_k = \int_{t_{k-1}
}^{t_k } {\left( {\left\| {\nabla \tilde {\omega }_1 }
\right\|_{L^2({\mathbb R^3} )}^2 +\left\| {\nabla \tilde {\omega }_2 }
\right\|_{L^2({\mathbb R^3} )}^2 +\left\| {\nabla \tilde {\omega }_3 }
\right\|_{L^2({\mathbb R^3} )}^2 } \right)}  \\
 &\qquad \qquad \qquad \qquad \qquad \qquad \qquad   k=1,\cdots ,N \\
\end{split}
\end{equation*}
then we have
\[
{\left( 1-4K_0 \frac{\Delta t_{k} }{\mu_\varepsilon} \right)}\;M_k\; + \;\delta_k \le\; M_{k-1}
\]

The partition is assumed to be fine enough. Because of the local existence of Galerkin solution in section 3 and the absolute continuity of integration with respect to $t$, it is valid that $\delta_k \to 0$
as $\Delta t_k \to 0$.
\\

We may first consider the case that
\[
{M_{k-1}} \frac{\Delta t_k }{\delta_k} \to 0, \quad \mbox{as} \;\;\Delta t_k \to 0
\]
which may be a subsequence $k'$, still denoted $k$. At this time, we can choose $\varepsilon_k$ on each $(t_{k-1},\;t_k)$ such that
\[
\mu_{\varepsilon_k} = 4K_0\,M_{k-1}\;\frac{\Delta t_k }{\delta_k} \quad \mbox{and} \quad 1-4K_0\,\frac{\Delta t_k }{\mu_{\varepsilon_k}} \ge \frac{1}{2}
\]
$\varepsilon = \mathop{\max }\limits_k \{\varepsilon_k \}$.
\\

Then we obtain
\[
{\left( 1-4K_0 \frac{\Delta t_k }{\mu_{\varepsilon_k}} \right)} \;M_k\; + \;\delta_k \;= {\left( 1- \frac{\delta_k }{M_{k-1}} \right)} \;M_k\; + \;\delta_k\; \le\; M_{k-1}
\]
it follows that $M_k \le M_{k-1}$.
\\

Otherwise, $\delta_k \le O(\Delta t_k) {M_{k-1}}$. In this case, a convolution inequality in [1] is applied to get
\begin{equation*}
\begin{split}
 &\left\| {\overline {\overline \omega}_i^k } \right\|_{L^2({\mathbb R^3} )}^2 \le \left\|
{J_\varepsilon * \left| {\overline \omega_i^k } \right|(x)}
\right\|_{L^2({\mathbb R}^3)}^2 \le \left\| {J_\varepsilon }
\right\|_{L^1({\mathbb R}^3)}^2 \left\|
{\overline \omega_i^k } \right\|_{L^2({\mathbb R^3} )}^2 \\
 &\quad \quad \quad \;\; = \left\| {\overline \omega_i^k } \right\|_{L^2({\mathbb R^3}
)}^2 \le \mathop {\sup }\limits_{(t_{k - 1} ,t_k )} \left\| {\tilde
\omega_i } \right\|_{L^2({\mathbb R^3} )}^2 \\
 \end{split}
\end{equation*}
This leads to that
\begin{equation*}
\begin{split}
 &\{\;\left\| {\nabla {\overline { \overline \omega }}_1^k }
\right\|_{L^2({\mathbb R^3} )}^2 + \left\| {\nabla {\overline
{\overline \omega }}_2^k } \right\|_{L^2({\mathbb R^3} )}^2 + \left\| {\nabla {\overline {\overline \omega }}_3^k }
\right\|_{L^2({\mathbb R^3} )}^2 \}  \\
 &\qquad  \le \frac{1}{\Delta t_k} \int_{t_{k-1}
}^{t_k } {\left( {\left\| {\nabla \tilde {\omega }_1 }
\right\|_{L^2({\mathbb R^3} )}^2 +\left\| {\nabla \tilde {\omega }_2 }
\right\|_{L^2({\mathbb R^3} )}^2 +\left\| {\nabla \tilde {\omega }_3 }
\right\|_{L^2({\mathbb R^3} )}^2 } \right)} \; \le O(1) {M_{k-1}} \\
\end{split}
\end{equation*}

Since these $(t_{k-1},\;t_k)$ are of finite length, the number of them is finite. According to Cauchy-Schwartz inequality, similar to (8), we have
\begin{equation*}
\begin{split}
 &\partial _t \int_{\mathbb R^3} {(\tilde {\omega }_1^2 + \tilde {\omega }_2^2 +
\tilde {\omega }_3^2 )} \;\; + \int_{\mathbb R^3} {[\,(\partial _{x_1 }
\tilde {\omega }_1 )^2 + (\partial _{x_2 } \tilde {\omega }_1 )^2 +
(\partial _{x_3
} \tilde {\omega }_1 )^2} \\
 &\quad \quad \quad \quad \quad \quad \quad \qquad \qquad  + (\partial
_{x_1 } \tilde {\omega }_2 )^2 + (\partial _{x_2 } \tilde {\omega
}_2 )^2 +
(\partial _{x_3 } \tilde {\omega }_2 )^2 \\
 &\quad \quad \quad \quad \quad \quad \quad \qquad \qquad  + (\partial
_{x_1 } \tilde {\omega }_3 )^2 + (\partial _{x_2 } \tilde {\omega
}_3 )^2 +
(\partial _{x_3 } \tilde {\omega }_3 )^2\,] \\
 &\le 4\,\{\,(\int_{\mathbb R^3} { {\overline {u}}_1^{k^4} } )^{\frac{1}{2}}(\int_{\mathbb R^3} {
{\overline {\overline \omega }}_1^{k^4} } )^{\frac{1}{2}} + (\int_{\mathbb R^3}
{ {\overline {u}}_1^{k^4} } )^{\frac{1}{2}}(\int_{\mathbb R^3} {
{\overline {\overline \omega }}_2^{k^4} } )^{\frac{1}{2}} + (\int_{\mathbb R^3}
{ {\overline {u}}_1^{k^4} } )^{\frac{1}{2}}(\int_{\mathbb R^3} {
{\overline {\overline \omega }}_3^{k^4} } )^{\frac{1}{2}}
\\
 &\qquad \;\;\; (\int_{\mathbb R^3} { {\overline {u}}_2^{k^4} } )^{\frac{1}{2}}(\int_{\mathbb R^3}
{ {\overline {\overline \omega }}_1^{k^4} } )^{\frac{1}{2}} +
(\int_{\mathbb R^3} { {\overline {u}}_2^{k^4} } )^{\frac{1}{2}}(\int_{\mathbb R^3}
{ {\overline {\overline \omega }}_2^{k^4} } )^{\frac{1}{2}} +
(\int_{\mathbb R^3} { {\overline {u}}_2^{k^4} } )^{\frac{1}{2}}(\int_{\mathbb R^3}
{ {\overline {\overline \omega }}_3^{k^4} } )^{\frac{1}{2}}
\\
 &\qquad \;\;\; (\int_{\mathbb R^3} { {\overline {u}}_3^{k^4} } )^{\frac{1}{2}}(\int_{\mathbb R^3}
{ {\overline {\overline \omega }}_1^{k^4} } )^{\frac{1}{2}} +
(\int_{\mathbb R^3} { {\overline {u}}_3^{k^4} } )^{\frac{1}{2}}(\int_{\mathbb R^3}
{ {\overline {\overline \omega }}_2^{k^4} } )^{\frac{1}{2}} +
(\int_{\mathbb R^3} { {\overline {u}}_3^{k^4} } )^{\frac{1}{2}}(\int_{\mathbb R^3}
{ {\overline {\overline \omega }}_3^{k^4} }
)^{\frac{1}{2}}\,\} \\
 &= 4\,\{\,\,\left\| { {\overline {u}}_1^k } \right\|_{L^4({\mathbb R^3} )}^2 (\,\left\|
{ {\overline {\overline \omega }}_1^k } \right\|_{L^4({\mathbb R^3} )}^2 +
\left\| { {\overline {\overline \omega }}_2^k } \right\|_{L^4({\mathbb R^3}
)}^2 + \left\| { {\overline {\overline \omega
}}_3^k } \right\|_{L^4({\mathbb R^3} )}^2 ) \\
 &\quad \;\; + \left\| { {\overline {u}}_2^k } \right\|_{L^4({\mathbb R^3} )}^2 (\,\left\|
{ {\overline {\overline \omega }}_1^k } \right\|_{L^4({\mathbb R^3} )}^2 +
\left\| { {\overline {\overline \omega }}_2^k } \right\|_{L^4({\mathbb R^3}
)}^2 + \left\| { {\overline {\overline \omega
}}_3^k } \right\|_{L^4({\mathbb R^3} )}^2 ) \\
 &\quad \;\; + \left\| { {\overline {u}
 }_3^k } \right\|_{L^4({\mathbb R^3} )}^2 (\,\left\|
{ {\overline {\overline \omega }}_1^k } \right\|_{L^4({\mathbb R^3} )}^2 +
\left\| { {\overline {\overline \omega }}_2^k } \right\|_{L^4({\mathbb R^3}
)}^2 + \left\| { {\overline {\overline \omega
}}_3^k } \right\|_{L^4({\mathbb R^3} )}^2 )\,\} \\
 &= 4\,(\;\left\| { {\overline {u}}_1^k } \right\|_{L^4({\mathbb R^3} )}^2 + \left\| { {\overline
{u}}_2^k } \right\|_{L^4({\mathbb R^3} )}^2 + \left\| { {\overline {u}}_3^k }
\right\|_{L^4({\mathbb R^3} )}^2 )\; (\;\left\| { {\overline
{\overline \omega }}_1^k } \right\|_{L^4({\mathbb R^3} )}^2 + \left\| {
{\overline {\overline \omega }}_2^k } \right\|_{L^4({\mathbb R^3} )}^2 + \left\|
{ {\overline {\overline \omega }}_3^k }
\right\|_{L^4({\mathbb R^3} )}^2 ) \\
 \end{split}
\end{equation*}
\\

From Sobolev imbedding theorem in [1], there exists a constant $C_1
> 0$ independent of $\omega $ such that
\[
\left\| { {\overline {\overline \omega }}_i^k } \right\|_{L^4({\mathbb R^3}
)}^2 \le C_1 \{\;\left\| { {\overline {\overline \omega }}_i^k }
\right\|_{L^2({\mathbb R^3} )}^2 + \left\| {\nabla  {\overline
{\overline \omega }}_i^k } \right\|_{L^2({\mathbb R^3} )}^2 \},\quad \quad i = 1,2,3
\]
Therefore,
\begin{equation*}
\begin{split}
 &\int_{\mathbb R^3} {(\tilde {\omega }_1^2 + \tilde {\omega }_2^2 + \tilde {\omega
}_3^2 )} \;\; + \int_{t_{k - 1} }^t {(\;\left\| {\nabla \tilde
{\omega }_1 } \right\|_{L^2({\mathbb R^3} )}^2 + \left\| {\nabla \tilde
{\omega }_2 } \right\|_{L^2({\mathbb R^3} )}^2 + \left\| {\nabla \tilde
{\omega }_3 }
\right\|_{L^2({\mathbb R^3} )}^2 )} \\
 &\le \int_{\mathbb R^3} {(\tilde {\omega }_1^{k - 1^2} + \tilde {\omega }_2^{k -
1^2} + \tilde {\omega }_3^{k - 1^2} )} \;\; + \\
 &+ \,C_2\,(\;\left\| {{\overline {u}}_1^k } \right\|_{L^2({\mathbb R^3} )}^2 + \left\|
{{\overline {u}}_2^k } \right\|_{L^2({\mathbb R^3} )}^2 + \left\| {{\overline
{u}}_3^k } \right\|_{L^2({\mathbb R^3} )}^2 + \left\| {\nabla {\overline
{u}}_1^k } \right\|_{L^2({\mathbb R^3} )}^2 + \left\| {\nabla {\overline
{u}}_2^k } \right\|_{L^2({\mathbb R^3} )}^2 + \left\| {\nabla {\overline
{u}}_3^k }
\right\|_{L^2({\mathbb R^3} )}^2 ) \\
 &\quad \quad \times \int_{t_{k - 1} }^t {\int_{\mathbb R^3} {( {\overline {\overline \omega
}}_1^{k^2} +  {\overline {\overline \omega }}_2^{k^2} +
{\overline {\overline \omega }}_3^{k^2}
)} } \;\; + \\
 &+ \,C_2\,\int_{t_{k - 1} }^t {(\;\left\| {{\overline {u}}_1^k }
\right\|_{L^2({\mathbb R^3} )}^2 + \left\| {{\overline {u}}_2^k }
\right\|_{L^2({\mathbb R^3} )}^2 + } \left\| {{\overline {u}}_3^k }
\right\|_{L^2({\mathbb R^3} )}^2 + \left\| {\nabla {\overline {u}}_1^k }
\right\|_{L^2({\mathbb R^3} )}^2 + \left\| {\nabla {\overline {u}}_2^k }
\right\|_{L^2({\mathbb R^3} )}^2 + \left\| {\nabla {\overline {u}}_3^k }
\right\|_{L^2({\mathbb R^3} )}^2 ) \\
 &\quad \quad \times \,(\;\left\| {\nabla  {\overline {\overline \omega }}_1^k }
\right\|_{L^2({\mathbb R^3} )}^2 + \left\| {\nabla  { \overline
{\overline \omega }}_2^k } \right\|_{L^2({\mathbb R^3} )}^2 + \left\| {\nabla
 {\overline {\overline \omega }}_3^k }
\right\|_{L^2({\mathbb R^3} )}^2 ) \\
 \end{split}
\end{equation*}
\\
Thus we have
\begin{equation*}
\begin{split}
 &\int_{\mathbb R^3} {(\tilde {\omega
}_1^2 + \tilde {\omega }_2^2 + \tilde {\omega }_3^2 )} \;\; +
\int_{t_{k - 1} }^{t } {(\;\left\| {\nabla \tilde {\omega }_1 }
\right\|_{L^2({\mathbb R^3} )}^2 + \left\| {\nabla \tilde {\omega }_2 }
\right\|_{L^2({\mathbb R^3} )}^2 +
\left\| {\nabla \tilde {\omega }_3 } \right\|_{L^2({\mathbb R^3} )}^2 )} \\
 &\qquad \le \int_{\mathbb R^3} {(\tilde {\omega }_1^{k - 1^2} + \tilde {\omega
}_2^{k - 1^2} + \tilde {\omega }_3^{k - 1^2} )} \;\; + \\
 &\qquad \quad + \,C_2\,\left( \;{\mathop {\sup }\limits_{(t_{k - 1} ,t_k )} \{\,\left\| {u_1
} \right\|_{L^2({\mathbb R^3} )}^2 + \left\| {u_2 } \right\|_{L^2({\mathbb R^3}
)}^2 +
\left\| {u_3 } \right\|_{L^2({\mathbb R^3} )}^2 \} \;\;+ } \right. \\
 &\quad \quad \quad \quad \quad \quad \quad  + \left.
{\frac{1}{\Delta t_k }\int_{t_{k - 1} }^{t_k } {\{\,\left\| {\nabla
u_1 } \right\|_{L^2({\mathbb R^3} )}^2 + \left\| {\nabla u_2 }
\right\|_{L^2({\mathbb R^3} )}^2
+ \left\| {\nabla u_3 } \right\|_{L^2({\mathbb R^3} )}^2 \}} } \right) \\
 &\quad \qquad \times \int_{t_{k - 1} }^{t } { \int_{\mathbb R^3} {(\tilde {\omega
}_1^2 + \tilde {\omega }_2^2 + \tilde {\omega }_3^2 )} } \;\; + \\
 &\quad \qquad  + \,C_3\,\left( {\Delta t_k \mathop {\sup
}\limits_{(t_{k - 1} ,t_k )} \{\,\left\| {u_1 } \right\|_{L^2({\mathbb R^3}
)}^2 + \left\| {u_2 } \right\|_{L^2({\mathbb R^3} )}^2 + \left\| {u_3 }
\right\|_{L^2({\mathbb R^3} )}^2 \} \;\;+ } \right. \\
 &\quad \quad \quad \quad \quad \quad \quad  + \left.
{\int_{t_{k - 1} }^{t_k } {\{\,\left\| {\nabla u_1 }
\right\|_{L^2({\mathbb R^3} )}^2 + \left\| {\nabla u_2 }
\right\|_{L^2({\mathbb R^3} )}^2 + \left\| {\nabla u_3
} \right\|_{L^2({\mathbb R^3} )}^2 \}} } \right) {M_{k-1}}\\
 \end{split}
\end{equation*}
where $C_2, C_3 >0$ are constants indepentent of $k$. Set
\begin{equation*}
\begin{split}
 &K_k^ * = \Delta t_k \mathop {\sup }\limits_{(t_{k - 1} ,t_k )} \{\,\left\|
{u_1 } \right\|_{L^2({\mathbb R^3} )}^2 + \left\| {u_2 }
\right\|_{L^2({\mathbb R^3} )}^2
+ \left\| {u_3 } \right\|_{L^2({\mathbb R^3} )}^2 \} + \\
 &\quad \quad \quad \quad \quad + \int_{t_{k - 1} }^{t_k } {\{\,\left\|
{\nabla u_1 } \right\|_{L^2({\mathbb R^3} )}^2 + \left\| {\nabla u_2 }
\right\|_{L^2({\mathbb R^3} )}^2 + \left\| {\nabla u_3 }
\right\|_{L^2({\mathbb R^3} )}^2
\}} \\
 \end{split}
\end{equation*}
and
\[
 f_k (t) = \,\mathop {\sup }\limits_{(t_{k - 1} ,t)}
\int_{\mathbb R^3} {(\tilde {\omega }_1^2 + \tilde {\omega }_2^2 + \tilde
{\omega
}_3^2 )}
\]
\\
Then we arrive at
\[
f_k (t) \le M_{k-1} \;+ C_2 \frac{1}{\Delta t_k } K_k^ *
\;\int_{t_{k-1}}^{t} {f_k (t)}  + \; C_3 K_k^ * {M_{k-1}}
\]
\\
By using Gronwall inequality it follows that
\[
M_k \le \left( 1 + C_3 K_k^ * \right)\; \exp \left( C_2 K_k^ *
\right)\,{M_{k-1}}
\]

Note that
\begin{equation*}
\begin{split}
 &\sum\limits_{k = 1}^N {K_k^ * } \le T\,\mathop {\sup }\limits_{t \in (0,T)}
\int_{\mathbb R^3} {(u_1^2 + u_2^2 + u_3^2 )} + \int_0^T {(\,\left\|
{\nabla u_1 } \right\|_{L^2({\mathbb R^3} )}^2 + \left\| {\nabla u_2 }
\right\|_{L^2({\mathbb R^3} )}^2
+ \left\| {\nabla u_3 } \right\|_{L^2({\mathbb R^3} )}^2 )} \\
 &\quad \qquad  \le (T+1)\,K_0 < + \infty \\
 \end{split}
\end{equation*}
\\
Hence, combining above two cases, we obtain
\\
\begin{equation*}
\begin{split}
 & M_1 \le \left( 1 + C_3 K_1^ * \right) \;\exp
\left( C_2 K_1^ * \right) \,M_0 \\
 & M_2 \le \left( 1 + C_3 K_1^ * \right)\,\left( 1 + C_3 K_2^ * \right) \;\exp
\left( C_2 \sum\limits_{k = 1}^2 {K_k^ * } \right) \,M_0 \\
 & \; \cdots \; \cdots \cdots \cdots \cdots \cdots \cdots \\
 & M_N \le \prod\limits_{k = 1}^N {(1 + C_3 K_k^\ast )} \;\exp
\left( C_2 \sum\limits_{k = 1}^N {K_k^ * } \right) \,M_0 \\
 \end{split}
\end{equation*}
Note that
\begin{equation*}
\begin{split}
 &\quad \prod\limits_{k = 1}^N {(1 + C_3 K_k^\ast )} = \exp \left( {\ln
\prod\limits_{k = 1}^N {(1 + C_3 K_k^\ast )} } \right)   \\
 &= \exp \left(
{\sum\limits_{k = 1}^N {\ln (1 + C_3 K_k^\ast )} } \right) \le \exp \left(
{C_3 \sum\limits_{k = 1}^N {K_k^\ast } } \right) = \exp (C_3 (T + 1)K_0 )   \\
 \end{split}
\end{equation*}
These mean that
\[
 M_k \le  M_0 \;\exp \left( (C_2 + C_3)\, (T+1)K_0 \right) \qquad   k=1,\cdots ,N
\]

Finally we get
\begin{equation*}
\begin{split}
&\mathop {\sup }\limits_{t \in (0,T)} \int_{\mathbb R^3} {(\tilde {\omega
}_1^2 + \tilde {\omega }_2^2 + \tilde {\omega }_3^2 )} \le
\, \mathop {\max }\limits_k \{M_k \} \\
&\qquad \qquad  \le  M_0 \;\exp \left( (C_2 + C_3) \,(T+1)K_0 \right)   \\
 \end{split}
\end{equation*}
\\

From condition (3) it is found that $M_0 $ is bounded, namely,
\begin{equation*}
\begin{split}
 &\int_{{\mathbb R}^3} {( {\omega }_{10}^2 (x)\, + {\omega }_{20}^2 (x)\,
+ {\omega }_{30}^2 (x))} \\
 &\quad \le \int\limits_{\left| x \right| \le R} {( {\omega }_{10}^2
(x)\, + {\omega }_{20}^2 (x)\, + {\omega }_{30}^2 (x))} \; +
\;12\,C_\mu ^2 \int\limits_{\left| x \right| > R} {\frac{1}{(1 +
\left| x
\right|)^{2\sigma }}} < + \infty \\
 \end{split}
\end{equation*}
and $R$ is a large constant.

This conclusion is also true for the weak solution of problem (7),
by means of the result of section 3 and the lower limit of Galerkin
sequence according to the page 196 of [4].
\\
\\
\\

\textbf{3. Existence}

In this section we have to consider the existence of solutions of
the auxiliary problems. We just need considering the following
system on $(0,\delta )$:
\begin{equation}
\label{eq9}
\begin{split}
 &\partial _t \tilde\omega _1 + \overline {u}_1 \partial _{x_1 } \overline{\overline {\omega }}_1 + \overline
{u}_2 \partial _{x_2 } \overline{\overline {\omega }}_1 + \overline
{u}_3
\partial _{x_3 } \overline{\overline {\omega }}_1 - \overline{\overline {\omega }}_1 \partial _{x_1 }
\overline {u}_1 - \overline{\overline {\omega }}_2 \partial _{x_2 }
\overline {u}_1 - \overline{\overline {\omega }}_3
\partial _{x_3 } \overline
{u}_1 +\partial _{x_1 } q = \Delta \tilde\omega _1 \\
 &\partial _t \tilde\omega _2 + \overline {u}_1 \partial _{x_1 } \overline{\overline {\omega }}_2 + \overline
{u}_2 \partial _{x_2 } \overline{\overline {\omega }}_2 + \overline
{u}_3
\partial _{x_3 } \overline{\overline {\omega }}_2 - \overline{\overline {\omega }}_1 \partial _{x_1 }
\overline {u}_2 - \overline{\overline {\omega }}_2
\partial _{x_2 } \overline {u}_2 - \overline{\overline {\omega }}_3 \partial _{x_3 } \overline {u}_2
+ \partial _{x_2 } q = \Delta \tilde\omega _2 \\
 &\partial _t \tilde\omega _3 + \overline {u}_1 \partial _{x_1 } \overline{\overline {\omega }}_3 + \overline
{u}_2 \partial _{x_2 } \overline{\overline {\omega }}_3 + \overline
{u}_3
\partial _{x_3 } \overline{\overline {\omega }}_3 - \overline{\overline {\omega }}_1 \partial
_{x_1 } \overline {u}_3 - \overline{\overline {\omega }}_2
\partial _{x_2 } \overline {u}_3 - \overline{\overline {\omega }}_3 \partial _{x_3 } \overline {u}_3
+\partial _{x_3 } q = \Delta \tilde\omega _3 \\
\end{split}
\end{equation}
with the initial value $ \tilde\omega _i (x,0)=\omega _{i0} \;\;(i=1,2,3)$ and
\[
\overline {u}_i (x)=\frac{1}{\delta }\int_0^\delta {u_i (x,t)dt}
\]
and
\[
\overline \omega_i (x)=\frac{1}{\delta }\int_0^\delta \tilde\omega_i
(x,t)dt,\quad \overline{\overline \omega}_i (x) = J_\varepsilon *
\overline \omega _i (x), \quad i=1,2,3
\]
as well as the incompressible conditions:
\begin{equation*}
\begin{split}
 &\partial _{x_1 } u_1 +\partial _{x_2 } u_2 + \partial _{x_3 } u_3
=0\quad \, \Rightarrow \quad \partial _{x_1 } \overline {u}_1
+\partial _{x_2 } \overline
{u}_2 + \partial _{x_3 } \overline {u}_3 = 0 \\
 &\partial _{x_1 } \tilde\omega _1 +\partial _{x_2 } \tilde\omega _2 +\partial _{x_3 }
\tilde\omega _3 = 0 \quad \Rightarrow \quad \partial _{x_1 }
\overline{\overline {\omega} }_1
+ \partial _{x_2 } \overline{\overline {\omega }}_2 + \partial _{x_3 } \overline{\overline {\omega }}_3 = 0 \\
\end{split}
\end{equation*}
\\

(i) The Galerkin procedure is applied. For each $m$ and $i=1,2,3$ we
define an approximate solution $(\tilde\omega _{1m} ,\;\tilde\omega _{2m}
,\;\tilde\omega _{3m} )$ as follows :
\[
\tilde\omega _{im} =\sum\limits_{j=1}^m {g_{ij} (t)w_{ij} }
\]
where $\{w_{i1} ,\;\cdots ,\;w_{im} ,\cdots \}$ is the basis of $W$,
and $W$= the closure of ${\cal V}$ in the Sobolev space
$W^{2,4}({\mathbb R}^3)$, which is separable and is dense in $V$.
Thus by means of weighted function $\theta _{r} $ introduced in
Section 1,
\begin{equation}
\label{eq10}
\begin{split}
 &(\theta _{r} \partial _t \tilde\omega _{im} ,\;w_{il} )+(\,\theta
_{r} \nabla \tilde\omega _{im} ,\;\nabla w_{il} )+(\,\nabla
\tilde\omega _{im} ,\;w_{il} \nabla \theta _{r} )\;\; + \\
 &+ (\theta
_{r} (\overline {u}\cdot \nabla )\overline{\overline {\omega }}_{im}
,\;w_{il} )-(\theta _{r} (\overline{\overline {\omega }}_m \cdot
\nabla ) \overline {u}_i
,\;w_{il} )=0  \\
\end{split}
\end{equation}
let $ r\to +\infty $ we get
\begin{equation}
\label{eq11}
\begin{split}
 &(\partial _t \tilde\omega _{im} ,\;w_{il} )+(\nabla \tilde\omega _{im} ,\;\nabla
w_{il} )+((\overline {u}\cdot \nabla ) \overline{\overline {\omega
}}_{im} ,\;w_{il}
)-((\overline{\overline {\omega }}_m \cdot \nabla ) \overline {u}_i ,\;w_{il} )= 0  \\
 &\qquad \qquad  t\in (0,\delta ),\quad \tilde\omega _{im} (0)=\omega _{i0}^m ,\quad l=1,\cdots ,m
 \\
\end{split}
\end{equation}
where $\omega _{i0}^m $ is the orthogonal projection in $H$ of
$\omega _{i0} $ onto the space spanned by $w_{i1} ,\;\cdots
,\;w_{im} $. Therefore,
\begin{equation*}
\begin{split}
 &\sum\limits_{j=1}^m {(w_{ij} ,\;w_{il} ){g}'_{ij} (t)} +\sum\limits_{j=1}^m
{(\nabla w_{ij} ,\;\nabla w_{il} )g_{ij} (t)} \;\; + \\
 &\quad \quad \quad \quad +\sum\limits_{j=1}^m {\{((\overline {u}(t)\cdot \nabla
)w_{ij}^ * ,\;w_{il} )-((w_j^ * \cdot \nabla )w_{il} ,\;\overline
{u}_i (t))\}} \;\overline
{g}_{ij} (t)= 0 \\
\end{split}
\end{equation*}
where $ w_j^ * = J_\varepsilon * w_j, \; w_{ij}^ * = J_\varepsilon
* w_{ij},\;\; \overline {g}_{ij} (t)= \frac{1}{\delta } \int_0^\delta {
g_{ij} (t)dt}$ and $u_i \in L^\infty (0,T;H)$ from Section 1 which
are determined by equations (1). Inverting the nonsingular matrix
with elements $(w_{ij} ,\;w_{il} ),\;\;1\le j,l\le m$, we can write
above system in the following form
\begin{equation}
\label{eq12}
\begin{split}
{g}'_{ij} (t)+ \sum\limits_{l=1}^m {\alpha _{ijl} \;g_{il} (t)}
+\sum\limits_{l=1}^m {\beta _{ijl} \;\overline {g}_{il} (t)} = 0
\end{split}
\end{equation}
where $\alpha _{ijl} ,\;\,\beta _{ijl} $ are constants.

The initial conditions are equivalent to
\[
g_{ij} (0)=g_{ij}^0
=\mbox{the}\;j^{\,th}\;\mbox{component}\;\mbox{of}\;\omega _{i0}^m
\]

We construct a sequence $\{g_{ij}^k \}$ by using a successive approximation :
\begin{equation*}
\begin{split}
 &{g_{ij}^1}^\prime = -\sum\limits_{l=1}^m {\alpha _{ijl} g_{il}^0 }
-\sum\limits_{l=1}^m {\beta _{ijl} \overline {g}_{il}^0 } \quad
\Rightarrow \quad g_{ij}^1 =g_{ij}^0 -\int_0^t {\left(
{\sum\limits_{l=1}^m {\alpha _{ijl} g_{il}^0 } + \sum\limits_{l=1}^m
{\beta _{ijl} \overline {g}_{il}^0 } } \right)}
\\
 &{g_{ij}^2}^\prime =-\sum\limits_{l=1}^m {\alpha _{ijl} g_{il}^1 }
-\sum\limits_{l=1}^m {\beta _{ijl} \overline {g}_{il}^1 } \quad
\Rightarrow \quad g_{ij}^2 =g_{ij}^0 -\int_0^t {\left(
{\sum\limits_{l=1}^m {\alpha _{ijl} g_{il}^1 } +\sum\limits_{l=1}^m
{\beta _{ijl} \overline {g}_{il}^1 } } \right)}
\\
 &\quad \quad \quad \cdots \cdots \cdots \cdots \\
 &{g_{ij}^k}^\prime =-\sum\limits_{l=1}^m {\alpha _{ijl} g_{il}^{k-1} }
-\sum\limits_{l=1}^m {\beta _{ijl} \overline {g}_{il}^{k-1} } \quad
\Rightarrow \quad g_{ij}^k =g_{ij}^0 -\int_0^t {\left(
{\sum\limits_{l=1}^m {\alpha _{ijl} g_{il}^{k-1} }
+\sum\limits_{l=1}^m {\beta _{ijl} \overline {g}_{il}^{k-1}
} } \right)} \\
\end{split}
\end{equation*}
so that
\[
\left| {g_{ij}^k (t)-g_{ij}^{k-1} (t)} \right|\le \int_0^t {\left(
{\sum\limits_{l=1}^m {\left| {\alpha _{ijl} } \right|\;\left|
{g_{il}^{k-1} (t)-g_{il}^{k-2} (t)} \right|} +\sum\limits_{l=1}^m
{\left| {\beta _{ijl} } \right|\;\left| {\overline {g}_{il}^{k-1}
(t)- \overline {g}_{il}^{k-2} (t)} \right|} } \right)}
\]

It follows that
\[
\mathop {\max }\limits_{i,j} \;\mathop {\sup }\limits_t \left| {g_{ij}^k
(t)-g_{ij}^{k-1} (t)} \right|\le \mathop {\max }\limits_{i,j}
\sum\limits_{l=1}^m {\left( {\left| {\alpha _{ijl} } \right|+\left| {\beta
_{ijl} } \right|} \right)\cdot t\cdot \mathop {\max }\limits_{i,j} \;\mathop
{\sup }\limits_t \;\left| {g_{ij}^{k-1} (t)-g_{ij}^{k-2} (t)} \right|}
\]

Taking $ \delta :\,=\frac{1}{\mathop {\max }\limits_{i,j}
\sum\limits_{l=1}^m {\left( {\left| {\alpha _{ijl} } \right|+2\left|
{\beta _{ijl} } \right|} \right)} } $, as $t\le \delta $, then
choosing $\delta ^\ast $:
\[
0<\delta ^\ast =\frac{\mathop {\max }\limits_{i,j} \sum\limits_{l=1}^m
{\left( {\left| {\alpha _{ijl} } \right|+\left| {\beta _{ijl} } \right|}
\right)} }{\mathop {\max }\limits_{i,j} \sum\limits_{l=1}^m {\left( {\left|
{\alpha _{ijl} } \right|+2\left| {\beta _{ijl} } \right|} \right)} }<1
\]
we have
\[
\mathop {\max }\limits_{i,j} \;\left\| {g_{ij}^k -g_{ij}^{k-1} }
\right\|_\infty \le \delta ^\ast \mathop {\max }\limits_{i,j} \;\left\|
{g_{ij}^{k-1} -g_{ij}^{k-2} } \right\|_\infty \le \cdots \le (\delta ^\ast
)^{k-1}\mathop {\max }\limits_{i,j} \;\left\| {g_{ij}^1 -g_{ij}^0 }
\right\|_\infty
\]

For any $n,\;k$ (we can set $n>k$ without loss of generality), we get
\begin{equation*}
\begin{split}
 &\mathop {\max }\limits_{i,j} \;\left\| {g_{ij}^n -g_{ij}^k }
\right\|_\infty \le \mathop {\max }\limits_{i,j} \;\left\| {g_{ij}^n
-g_{ij}^{n-1} } \right\|_\infty +\cdots +\mathop {\max }\limits_{i,j}
\;\left\| {g_{ij}^{k+1} -g_{ij}^k } \right\|_\infty \\
 &\le ((\delta ^\ast )^{n-1}+\cdots +(\delta ^\ast )^k)\;\,\mathop {\max
}\limits_{i,j} \;\left\| {g_{ij}^1 -g_{ij}^0 } \right\|_\infty =(\delta
^\ast )^k\frac{1-(\delta ^\ast )^{n-k}}{1-\delta ^\ast }\mathop {\max
}\limits_{i,j} \;\left\| {g_{ij}^1 -g_{ij}^0 } \right\|_\infty \\
 &\to 0\quad (k\to \infty ) \\
\end{split}
\end{equation*}
Thus, for every $i=1,2,3;\;\;j=1,\cdots ,m$, $\{g_{ij}^k \}$ is a
Cauchy sequence in $L^\infty (0,\delta )$. Since $L^\infty (0,\delta
)$ is complete, then there exists a function $g_{ij}^\ast \in
L^\infty (0,\delta )$ such that $\left\| {g_{ij}^k -g_{ij}^\ast }
\right\|_\infty \to 0$ as $k\to \infty $.
\\

From
\[
g_{ij}^k (t)=g_{ij}^0 -\int_0^t {\left( {\sum\limits_{l=1}^m {\alpha
_{ijl} \;g_{il}^{k-1} (t)} +\sum\limits_{l=1}^m {\beta _{ijl}
\;\overline {g}_{il}^{k-1} (t)} } \right)}
\]
let $k\to \infty $, it follows that
\[
g_{ij}^\ast (t)=g_{ij}^0 -\int_0^t {\left( {\sum\limits_{l=1}^m
{\alpha _{ijl} \;g_{il}^\ast (t)} +\sum\limits_{l=1}^m {\beta _{ijl}
\;\overline {g}_{il}^\ast (t)} } \right)}
\]
i.e., $g_{ij}^\ast $ is a solution of the system (12) on $(0,\delta )$ for
which $g_{ij}^\ast (0)=g_{ij}^0 $, $i=1,2,3;\;\;j=1,\cdots ,m$.
\\

(ii) By means of the weighted function $\theta _{r} $ :
\begin{equation*}
\begin{split}
 &\sum\limits_{i=1}^3 {(\,\theta _{r} \partial _t \tilde\omega _{im}
,\;\tilde\omega _{im} )} + \sum\limits_{i=1}^3 {(\,\theta _{r} \nabla
\tilde\omega _{im} ,\;\nabla \tilde\omega _{im} )} +\sum\limits_{i=1}^3
{(\,\nabla
\tilde\omega _{im} ,\;\,\tilde\omega _{im} \nabla \theta _{r} )} \;\; + \\
 &\quad \quad \;+\sum\limits_{i=1}^3 {(\,\theta _{r} (\overline {u}\cdot
\nabla )\overline{\overline {\omega }}_{im} ,\;\tilde\omega _{im} )}
-\sum\limits_{i=1}^3 {(\,\theta _{r} (\overline{\overline {\omega
}}_m \cdot \nabla ) \overline {u}_i
,\;\tilde\omega _{im} )} = 0 \\
\end{split}
\end{equation*}
Let $ r\to +\infty $ we get
\[
\sum\limits_{i=1}^3 {(\partial _t \tilde\omega _{im} ,\;\tilde\omega _{im} )}
+\sum\limits_{i=1}^3 {(\nabla \tilde\omega _{im} ,\;\nabla \tilde\omega _{im} )}
+\sum\limits_{i=1}^3 {((\overline {u}\cdot \nabla )
\overline{\overline {\omega }}_{im} ,\;\tilde\omega _{im} )}
-\sum\limits_{i=1}^3 {((\overline{\overline {\omega }}_m \cdot
\nabla ) \overline {u}_i ,\;\tilde\omega _{im} )} = 0
\]
Then we write
\begin{equation*}
\begin{split}
 &\frac{1}{2}\frac{d}{dt}\left( {\sum\limits_{i=1}^3 {\left\| {\tilde\omega _{im} }
\right\|_{L^2({\mathbb R}^3)}^2 } } \right) + \sum\limits_{i=1}^3
{\left\| {\nabla \tilde\omega _{im} } \right\|_{L^2({\mathbb R}^3)}^2 }
-\sum\limits_{i=1}^3 {((\overline
{u}\cdot \nabla )\tilde\omega _{im} ,\;\overline{\overline {\omega} }_{im} )} \;\; + \\
 &\quad \quad \quad \quad \quad + \sum\limits_{i=1}^3 {((\overline{\overline {\omega }}_m
\cdot \nabla )\tilde\omega _{im} ,\;\overline {u}_i )} = 0 \\
\end{split}
\end{equation*}

Similar to those in the section 2, and $\eta $ is
chosen to be small enough, we have
\[
\sum\limits_{i=1}^3 {\left\| {\tilde\omega _{im} } \right\|_{L^2({{\mathbb
R}^3} )}^2 } + \int_0^\eta {\left( {\sum\limits_{i=1}^3 {\left\|
{\nabla \tilde\omega _{im} } \right\|_{L^2({{\mathbb R}^3} )}^2 } }
\right)} \le 2\, \left(
{\sum\limits_{i=1}^3 {\left\| {\omega _{i0}^m }
\right\|_{L^2({{\mathbb R}^3} )}^2 } } \right)
\]
as $1-4\,K_0 \,\eta/\mu_\varepsilon \ge 1/2 $ . Hence,
\begin{equation}
\label{eq13} \mathop {\sup }\limits_{t\in (0,\eta )} \left(
{\sum\limits_{i=1}^3 {\left\| {\tilde\omega _{im} }
\right\|_{L^2({{\mathbb R}^3} )}^2 } } \right)\le
2\, \left( {\sum\limits_{i=1}^3 {\left\|
{\omega _{i0} } \right\|_{L^2({{\mathbb R}^3} )}^2 } } \right)
\end{equation}
and
\begin{equation}
\label{eq14} \sum\limits_{i=1}^3 {\left\| {\tilde\omega _{im} (\eta )}
\right\|_{L^2({{\mathbb R}^3} )}^2 } + \int_0^\eta {\left(
{\sum\limits_{i=1}^3 {\left\| {\nabla \tilde\omega _{im} }
\right\|_{L^2({{\mathbb R}^3} )}^2 } } \right)} \le
2\,\left( {\sum\limits_{i=1}^3 {\left\|
{\omega _{i0} } \right\|_{L^2({{\mathbb R}^3} )}^2 } } \right)
\end{equation}

The inequalities (13) and (14) are valid for any fixed $\delta \le
\eta $.
\\

(iii) Let $ {\mathop \omega \limits^\circ }_m $ denote the function from ${\mathbb
R}$ into $V$, which is equal to $\tilde\omega _m $ on $(0,\delta )$ and to
0 on the complement of this interval. The Fourier transform of
$ {\mathop \omega \limits^\circ }_m $ is denoted by $\hat {\omega }_m $. We want to
show that
\[
\int_{-\infty }^{+\infty } {\left| \tau \right|^{2\gamma }\left(
{\sum\limits_{i=1}^3 {\left\| {\hat {\omega }_{im} (\tau )}
\right\|_{L^2(\Omega )}^2 } } \right)} \,d\tau <+\infty ,\quad \quad
\forall\; \Omega \subset {\mathbb R}^3
\]
For some $\gamma >0$. Along with (14) this will imply that
\\

$ {\mathop \omega \limits^\circ }_m $ \textit{belongs to a bounded set of} $H^\gamma
({\mathbb R},\;H^1(\Omega ),\;L^2(\Omega )),\quad \forall\; \Omega $
\\
\\
and will enable us to apply the result of compactness.

We observe that (10) can be written as
\begin{equation*}
\begin{split}
 &\frac{d}{dt}\left( {\sum\limits_{i=1}^3 {(\,\theta _{r} {\mathop \omega \limits^\circ }_{im} ,\;w_{ij} )} } \right)=\sum\limits_{i=1}^3 {(\,\theta
_{r} {\mathop f \limits^\circ }_{im} ,\;w_{ij} )} +\sum\limits_{i=1}^3 {(\,\theta
_{r} \omega _{i0}^m ,\;w_{ij} )\,} \eta _0 \; - \\
 &\qquad \qquad \qquad \quad \quad \quad \quad \quad \quad -\sum\limits_{i=1}^3
{(\,\theta _{r} \tilde\omega _{im} (\delta ),\;w_{ij} )\,} \eta _\delta
\\
\end{split}
\end{equation*}
where $\eta _0 ,\;\eta _\delta $ are Dirac distributions at 0 and
$\delta $, and
\begin{equation*}
\begin{split}
 &f_{im} =-\Delta \tilde\omega _{im} +(\overline {u}\cdot \nabla ) \overline{\overline {\omega
}}_{im} -(\overline{\overline {\omega }}_m \cdot \nabla )\;\overline
{u}_i
 \\
 &{\mathop f \limits^\circ }_{im} = f_{im} \;\; \mbox{on}\;\; (0,\delta ),\quad 0\; \mbox{ outside
this interval}  \\
\end{split}
\end{equation*}

By the Fourier transform,
\begin{equation*}
\begin{split}
 &2\mbox{i}\pi \tau \sum\limits_{i=1}^3 {(\,\theta _{r} \hat {\omega
}_{im} ,\;w_{ij} )} =\sum\limits_{i=1}^3 {(\,\theta _{r} \hat
{f}_{im} ,\;w_{ij} )} +\sum\limits_{i=1}^3 {(\,\theta _{r} \omega
_{i0}^m ,\;w_{ij} )} \; - \\
 &\quad \quad \quad \quad \quad \quad \quad \quad \quad -\sum\limits_{i=1}^3
{(\,\theta _{r} \tilde\omega _{im} (\delta ),\;w_{ij} )\,} \exp
(-2\mbox{i}\pi \delta \tau ) \\
\end{split}
\end{equation*}
where $\hat {\omega }_{im} $ and $\hat {f}_{im} $ denote the Fourier
transforms of ${\mathop \omega \limits^\circ }_{im} $ and ${\mathop f \limits^\circ }_{im} $
respectively.
\\

We multiply above equalities by $\hat {g}_{ij} (\tau )= $ Fourier
transform of ${\mathop g \limits^\circ }_{ij} $ and add the resulting equations for
$j=1,\cdots ,m$, we get
\begin{equation*}
\begin{split}
 &2\mbox{i}\pi \tau \sum\limits_{i=1}^3 {\left\| {\,\theta^{1/2} _{
r} \;\hat {\omega }_{im} (\tau )} \right\|_{L^2({\mathbb R}^3)}^2 }
=\sum\limits_{i=1}^3 {(\,\theta _{r} \hat {f}_{im} (\tau
),\;\hat {\omega }_{im} (\tau ))}  \\
 &\quad \quad +\sum\limits_{i=1}^3 {(\,\theta _{r} \omega _{i0}^m
,\;\hat {\omega }_{im} (\tau ))} -\sum\limits_{i=1}^3 {(\,\theta
_{r} \tilde\omega _{im} (\delta ),\;\hat {\omega }_{im} (\tau ))\,} \exp
(-2\mbox{i}\pi \delta \tau )
 \\
\end{split}
\end{equation*}

For some $\varphi _i \in V$,
\begin{equation*}
\begin{split}
 &\int_0^\delta {\sum\limits_{i=1}^3 {(\,\theta _{r} f_{im}
,\;\varphi _i )} } = \int_0^\delta {\sum\limits_{i=1}^3 {(-\theta
_{r} \Delta \tilde\omega _{im} ,\;\varphi _i )} } + \int_0^\delta
{\sum\limits_{i=1}^3 {(\,\theta _{r} (\overline {u}\cdot \nabla
)\overline{\overline
{\omega }}_{im} ,\;\varphi _i )} } \; - \\
 &\quad \quad \quad \quad \quad \quad {\kern 1pt}\qquad \qquad - \int_0^\delta
{\sum\limits_{i=1}^3 {(\,\theta _{r} (\overline{\overline {\omega
}}_m \cdot \nabla
)\;\overline {u}_i ,\;\varphi _i )} } \\
 &\quad =\int_0^\delta {\sum\limits_{i=1}^3 {(\,\theta _{r} \nabla \tilde\omega
_{im} ,\;\nabla \varphi _i )} } + \int_0^\delta {\sum\limits_{i=1}^3
{(\nabla
\tilde\omega _{im} ,\;\,\varphi _i \nabla \theta _{r} )} } \\
 &\qquad -\int_0^\delta {\sum\limits_{i=1}^3 {(\,\theta _{r} (\overline
{u}\cdot \nabla )\varphi _i ,\;\overline{\overline {\omega }}_{im}
)} } - \int_0^\delta
{\sum\limits_{i=1}^3 {(\,\varphi _i (\overline {u}\cdot \nabla )\theta _{r} ,\;\overline{\overline {\omega }}_{im} )} } \\
 &\qquad +\int_0^\delta {\sum\limits_{i=1}^3 {(\,\theta _{r} (\overline{\overline
{\omega }}_m \cdot \nabla )\varphi _i ,\;\overline {u}_i )} } +
\int_0^\delta {\sum\limits_{i=1}^3 {(\,\varphi _i
(\overline{\overline {\omega }}_m \cdot \nabla )\,\theta
_{r} ,\;\overline {u}_i )} } \\
\end{split}
\end{equation*}
Let $ r\to +\infty $ we get
\begin{equation*}
\begin{split}
 &\int_0^\delta {\sum\limits_{i=1}^3 {(f_{im} ,\;\varphi _i )} }
=\int_0^\delta {\sum\limits_{i=1}^3 {(\nabla \tilde\omega _{im} ,\;\nabla
\varphi _i )} } -\int_0^\delta {\sum\limits_{i=1}^3 {((\overline
{u}\cdot \nabla )\varphi _i ,\;\overline{\overline {\omega }}_{im}
)} } + \int_0^\delta {\sum\limits_{i=1}^3 {((\overline{\overline
{\omega }}_m \cdot \nabla )\varphi _i ,\;\overline {u}_i )} } \\
 &\le \int_0^\delta {\sum\limits_{i=1}^3 {\left\| {\nabla \tilde\omega _{im} }
\right\|_{L^2({\mathbb R}^3)} \left\| {\nabla \varphi _i }
\right\|_{L^2({\mathbb
R}^3)} } } \;\; + \\
\end{split}
\end{equation*}
\begin{equation*}
\begin{split}
 &\quad + 2\, \int_0^\delta \left( {\sum\limits_{i=1}^3 {\left\| {\overline {u}_i }
\right\|_{L^4({\mathbb R}^3 )}^2 } } \right)^{1/2}\left(
{\sum\limits_{i=1}^3 {\left\| {\overline{\overline {\omega }}_{im} }
\right\|_{L^4({\mathbb R}^3 )}^2 } } \right)^{1/2}\left(
{\sum\limits_{i=1}^3 {\left\| {\nabla \varphi _i }
\right\|_{L^2({\mathbb R}^3 )}^2 } } \right)^{1/2} \\
 &\le \int_0^\delta {\left( {\sum\limits_{i=1}^3 {\left\| {\nabla \tilde\omega
_{im} } \right\|_{L^2({\mathbb R}^3)}^2 } } \right)^{1/2}} \left(
{\sum\limits_{i=1}^3 {\left\| {\nabla \varphi _i }
\right\|_{L^2({\mathbb
R}^3)}^2 } } \right)^{1/2}+ \\
 &\quad +2C \delta \left( { {\sum\limits_{i=1}^3
{\{\,\left\| {\overline {u}_i } \right\|_{L^2({\mathbb R}^3)}^2
+\left\| {\nabla \overline
{u}_i } \right\|_{L^2({\mathbb R}^3)}^2 \}} } } \right)^{1/2} \\
 &\quad \quad \times \left( { {\sum\limits_{i=1}^3 {\{\,\left\|
{\overline {\omega }_{im} } \right\|_{L^2({\mathbb R}^3)}^2 +
\left\| {\nabla \overline {\omega }_{im} } \right\|_{L^2({\mathbb
R}^3)}^2 \}} } } \right)^{1/2}\left( {\sum\limits_{i=1}^3 {\left\|
{\nabla \varphi _i } \right\|_{L^2({\mathbb
R}^3)}^2 } } \right)^{1/2} \\
 &\le \int_0^\delta {\left( {\sum\limits_{i=1}^3 {\left\| {\nabla \tilde\omega
_{im} } \right\|_{L^2({\mathbb R}^3)}^2 } } \right)^{1/2}} \left\|
{\nabla
\varphi } \right\|_V + \\
 &\quad + 2C \left( {\delta \;\mathop {\sup }\limits_{(0,\delta )}
\;\sum\limits_{i=1}^3 {\left\| {u_i } \right\|_{L^2({\mathbb
R}^3)}^2 } +\int_0^\delta {\sum\limits_{i=1}^3 {\left\| {\nabla u_i
}
\right\|_{L^2({\mathbb R}^3)}^2 } } } \right)^{1/2} \\
 &\qquad \times \left( {\delta \;\mathop {\sup }\limits_{(0,\delta )}
\sum\limits_{i=1}^3 {\left\| {\tilde\omega _{im} } \right\|_{L^2({\mathbb
R}^3)}^2 } +\int_0^\delta {\sum\limits_{i=1}^3 {\left\| {\nabla
\tilde\omega _{im} } \right\|_{L^2({\mathbb R}^3)}^2 } } }
\right)^{1/2}\left\| {\nabla \varphi }
\right\|_V \\
\end{split}
\end{equation*}
this remains bounded according to (6) and (13), (14). Therefore,
\[
\int_0^\delta {\left\| {f_{im} (t)} \right\|_V dt} = \int_0^\delta
{\;\mathop {\sup }\limits_{\left\| \varphi \right\|_V =1}
\;\sum\limits_{i=1}^3 {(f_{im} ,\;\varphi _i )} } <+\infty
\]
it follows that
\[
\mathop {\sup }\limits_{\tau \in {\mathbb R}} \left\| {\hat {f}_{im}
(\tau )} \right\|_V <+\infty ,\quad \;\forall m
\]

Due to (13) we have
\[
\left\| {\omega _{im} (0)} \right\|_{L^2({\mathbb R}^3)} <+\infty
,\quad \quad \left\| {\tilde\omega _{im} (\delta )} \right\|_{L^2({\mathbb
R}^3)} <+\infty
\]
then by Poincare inequality,
\begin{equation*}
\begin{split}
 &\left| \tau \right|\;\sum\limits_{i=1}^3 {\left\| {\,\theta^{1/2}
_{r} \;\hat {\omega }_{im} (\tau )} \right\|_{L^2({\mathbb R}^3)}^2
} \le c_1 \sum\limits_{i=1}^3 {\left\| {\hat {f}_{im} (\tau )}
\right\|_V \;\left\| {\,\theta
_{r} \hat {\omega }_{im} (\tau )} \right\|_V } \\
 & \qquad \qquad \qquad \qquad \qquad \qquad + c_2
\sum\limits_{i=1}^3 {\left\| {\,\theta _{r} \hat {\omega }_{im}
(\tau )}
\right\|_{L^2({\mathbb R}^3)} }  \\
 &\quad \le c_3 \sum\limits_{i=1}^3 {\left\| {\nabla (\theta _{r} \hat
{\omega }_{im} (\tau ))} \right\|_{L^2({\mathbb R}^3)} }  \\
 &\quad \le c_4 \sum\limits_{i=1}^3 {\left( {\left\| {\,\hat {\omega }_{im} \nabla
\theta _{r} } \right\|_{L^2({\mathbb R}^3)} } \right.} +\left.
{\left\| {\,\theta _{r} \nabla \hat {\omega }_{im} }
\right\|_{L^2({\mathbb R}^3)} } \right)  \\
\end{split}
\end{equation*}

Using $x^2e^{-\kappa x}\le C_1 \; (\kappa>0)$ and assuming that $ r
$ is sufficiently large, we get
\begin{equation}
\label{eq15}
\begin{split}
 &\left| \tau \right|\;\sum\limits_{i=1}^3 {\left\| {\,\theta^{1/2}
_{r} \;\hat {\omega }_{im} (\tau )}
\right\|_{L^2({\mathbb R}^3)}^2 } \\
 &\quad \le c_5 \sum\limits_{i=1}^3
{\left\| {\,\theta^{1/2} _{r} \hat {\omega }_{im} }
\right\|_{L^2({\mathbb R}^3)} } +\;c_6 \;\sum\limits_{i=1}^3
{\left\| {\,\theta _{r} \nabla \hat {\omega }_{im} }
\right\|_{L^2({\mathbb R}^3)} } \\
\end{split}
\end{equation}

For $\gamma $ fixed, $\gamma <1/4$, we observe that
\[
\left| \tau \right|^{2\gamma }\le c_7 (\gamma )\frac{1+\left| \tau
\right|}{1+\left| \tau \right|^{1-2\gamma }},\quad \quad \forall
\tau \in {\mathbb R}
\]

Thus by (15),
\begin{equation*}
\begin{split}
 &\int_{-\infty }^{+\infty } {\left| \tau \right|^{2\gamma }\left(
{\sum\limits_{i=1}^3 {\left\| {\,\theta^{1/2} _{r} \;\hat {\omega
}_{im} (\tau )} \right\|_{L^2({\mathbb R}^3)}^2 } } \right)} \,d\tau
\le c_7 (\gamma )\int_{-\infty }^{+\infty } {\frac{1+\left| \tau
\right|}{1+\left| \tau \right|^{1-2\gamma }}\left(
{\sum\limits_{i=1}^3 {\left\| {\,\theta^{1/2} _{r} \;\hat {\omega
}_{im} (\tau )} \right\|_{L^2({\mathbb R}^3)}^2 } } \right)}
\,d\tau  \\
 &\quad \le c_8 \;\int_{-\infty }^{+\infty } {\frac{1}{1+\left| \tau
\right|^{1-2\gamma }}\;\sum\limits_{i=1}^3 {\left\| {\,\theta^{1/2}
_{r} \hat {\omega }_{im} (\tau )}
\right\|_{L^2({\mathbb R}^3)} } } d\tau \;\; +  \\
 &\qquad +\;c_9 \int_{-\infty }^{+\infty } {\frac{1}{1+\left| \tau \right|^{1-2\gamma
}}\;\sum\limits_{i=1}^3 {\left\| {\,\theta _{r} \nabla \hat {\omega
}_{im} (\tau )} \right\|_{L^2({\mathbb R}^3)} } } d\tau +\;c_{10}
\int_{-\infty }^{+\infty } {\sum\limits_{i=1}^3 {\left\|
{\,\theta^{1/2} _{r} \;\hat {\omega }_{im} (\tau )}
\right\|_{L^2({\mathbb R}^3)}^2 } } \,d\tau  \\
\end{split}
\end{equation*}
Because of the Parseval equality,
\begin{equation*}
\begin{split}
 &\int_{-\infty }^{+\infty } {\sum\limits_{i=1}^3 {\left\| {\,\theta
_{r} \hat {\omega }_{im} (\tau )} \right\|_{L^2({\mathbb R}^3)}^2 }
} \,d\tau =\int_0^\delta {\sum\limits_{i=1}^3 {\left\| {\,\theta
_{r}
\tilde\omega _{im} (t)} \right\|_{L^2({\mathbb R}^3)}^2 } } \,dt \\
 &\quad \quad \quad \quad \quad \quad \quad \quad \quad \quad \;\,\le
 C_3
\delta \;\mathop {\sup }\limits_{(0,\delta )} \;\sum\limits_{i=1}^3 {\left\|
{\tilde\omega _{im} } \right\|_{L^2({\mathbb R}^3)}^2 } < + \infty \\
 &\int_{-\infty }^{+\infty } {\sum\limits_{i=1}^3 {\left\| {\,\theta
_{r} \nabla \hat {\omega }_{im} (\tau )} \right\|_{L^2({\mathbb
R}^3)}^2 } } \,d\tau = \int_0^\delta {\sum\limits_{i=1}^3 {\left\|
{\,\theta _{r} \nabla \tilde\omega _{im} (t)} \right\|_{L^2({\mathbb
R}^3)}^2 } } \,dt
\\
 &\quad \quad \quad \quad \quad \quad \quad \quad \quad
\quad \,\le C_4 \int_0^\delta {\sum\limits_{i=1}^3 {\left\| {\nabla
\tilde\omega
_{im} } \right\|_{L^2({\mathbb R}^3)}^2 } } <+\infty \\
\end{split}
\end{equation*}
as $m\to \infty $. By Cauchy-Schwarz inequality and the Parseval
equality,
\begin{equation*}
\begin{split}
 &\int_{-\infty }^{+\infty } {\frac{1}{1+\left| \tau \right|^{1-2\gamma
}}\sum\limits_{i=1}^3 {\left\| {\,\theta^{1/2} _{r} \;\hat
{\omega }_{im} (\tau )} \right\|_{L^2({\mathbb R}^3)} } } d\tau \\
 &\quad \quad \le \sqrt 3 \left( {\int_{-\infty }^{+\infty }
{\frac{1}{(1+\left| \tau \right|^{1-2\gamma })^2}d\tau } }
\right)^{1/2}\left( {\int_0^\delta {\sum\limits_{i=1}^3 {\left\|
{\,\theta^{1/2} _{r} \;\tilde\omega _{im} (t)} \right\|_{L^2({\mathbb
R}^3)}^2 } }
dt} \right)^{1/2}<+\infty \\
\end{split}
\end{equation*}
\begin{equation*}
\begin{split}
 &\int_{-\infty }^{+\infty } {\frac{1}{1+\left| \tau \right|^{1-2\gamma
}}\sum\limits_{i=1}^3 {\left\| {\,\theta _{r} \nabla \hat {\omega
}_{im} (\tau )} \right\|_{L^2({\mathbb R}^3)} } } d\tau \\
 &\quad \quad \le \sqrt 3 \left( {\int_{-\infty }^{+\infty }
{\frac{1}{(1+\left| \tau \right|^{1-2\gamma })^2}d\tau } }
\right)^{1/2}\left( {\int_0^\delta {\sum\limits_{i=1}^3 {\left\|
{\,\theta _{r} \nabla \tilde\omega _{im} (t)} \right\|_{L^2({\mathbb
R}^3)}^2 } } dt}
\right)^{1/2}<+\infty \\
\end{split}
\end{equation*}
as $m\to \infty $ by $\gamma <1/4$ and (14).
\\

(iv) The estimates (13) and (14) enable us to assert the existence
of an element $\tilde\omega ^\ast \in L^2(0,\delta ;H^1(\Omega ))\cap
L^\infty (0,\delta ;L^2(\Omega )),\quad \forall\; \Omega \subset
{\mathbb R}^3$, and a subsequence $\tilde\omega _{{m}'} $ such that
\\

$\tilde\omega _{{m}'} \to \tilde\omega ^\ast $ \textit{in} $L^2(0,\delta
;H^1(\Omega ))$ \textit{weakly}, \textit{and in} $L^\infty (0,\delta
;L^2(\Omega ))$ \textit{weak-star},

\textit{as} ${m}'\to \infty $, \textit{for any} $\Omega \subset
{\mathbb R}^3$
\\

Due to (iii) we also have
\\

$\tilde\omega _{{m}'} \to \tilde\omega ^\ast $ \textit{in} $L^2(0,\delta
;L^2(\Omega ))$ \textit{strongly as} ${m}'\to \infty $, \textit{for
any} $\Omega \subset {\mathbb
R}^3$
\\

\noindent  which means
\\

$\tilde\omega _{{m}'} \to \tilde\omega ^\ast $ \textit{in} $L^2(0,\delta ;L_{
\mbox{
\begin{footnotesize}loc \end{footnotesize}} } ^2 (\Omega ))$ \textit{strongly}
\\

\noindent In particular, for a fixed $j$
\\

$\left. {\tilde\omega _{{m}'} } \right|_{{\Omega }'} \to \left. {\tilde\omega
^\ast } \right|_{{\Omega }'} $ \textit{in} $L^2(0,\delta
;L^2({\Omega }'))$
\textit{strongly}
\\

\noindent where ${\Omega }'$ denotes the support of $w_{ij} $. This
convergence result enable us to pass to the limit.
\\

Let $\psi _i $ be a continuously differentiable function on $(0,\delta )$
with $\psi _i (\delta )=0$. We multiply (11) by $\psi _i (t)$ then integrate
by parts. This leads to the equation
\begin{equation*}
\begin{split}
 &-\int_0^\delta {\sum\limits_{i=1}^3 {(\tilde\omega _{im} (t),\;\partial _t \psi
_i (t)w_{ij} )\,dt} } + \int_0^\delta {\sum\limits_{i=1}^3 {(\nabla \tilde\omega
_{im} ,\;\psi _i (t)\nabla w_{ij} )\,dt} } \\
 &+ \int_0^\delta {\sum\limits_{i=1}^3 {((\overline {u}\cdot \nabla )\overline{\overline {\omega
}}_{im} ,\;w_{ij} \psi _i (t))} } - \int_0^\delta
{\sum\limits_{i=1}^3 {((\overline{\overline {\omega} }_m \cdot
\nabla )\,\overline {u}_i ,\;w_{ij} \psi _i (t))} }
=\sum\limits_{i=1}^3 {(\omega _{i0}^m ,\;w_{ij} )\psi _i (0)} \\
\end{split}
\end{equation*}

Since $\tilde\omega _{i{m}'} $ converges to $\tilde\omega _i^\ast $ in
$L^2(0,\delta ;L^2(\Omega ))$ strongly as ${m}'\to \infty $, then
$\overline{\overline {\omega }}_{i{m}'} $ also converges strongly to
$\overline{\overline {\omega}}_i^\ast $, and
\begin{equation*}
\begin{split}
 &\int_0^\delta {\sum\limits_{i=1}^3 {(\tilde\omega _{i{m}'} ,\;\partial _t
\psi _i (t)w_{ij} )\,dt} } \to \int_0^\delta {\sum\limits_{i=1}^3
{(\tilde\omega _i^\ast ,\;\partial _t \psi _i (t)w_{ij} )\,dt} }  \\
 &\int_0^\delta {\sum\limits_{i=1}^3 {(\nabla \tilde\omega _{i{m}'} ,\;\psi _i
(t)\nabla w_{ij} )\,dt} } = -\int_0^\delta {\sum\limits_{i=1}^3
{(\tilde\omega
_{i{m}'} ,\;\psi _i (t)\Delta w_{ij} )\,dt} } \\
\end{split}
\end{equation*}
\begin{equation*}
\begin{split}
 &\quad \;\quad \quad \quad \to -\int_0^\delta {\sum\limits_{i=1}^3 {(\tilde\omega
_i^\ast ,\;\psi _i (t)\Delta w_{ij} )} } = \int_0^\delta
{\sum\limits_{i=1}^3
{(\nabla \tilde\omega _i^\ast ,\;\psi _i (t)\nabla w_{ij} )\,dt} } \\
 &\int_0^\delta {\sum\limits_{i=1}^3 {((\overline {u}\cdot \nabla ) \overline{\overline {\omega
}}_{i{m}'} ,\;w_{ij} \psi _i (t))} } = -\int_0^\delta
{\sum\limits_{i=1}^3
{((\overline {u}\cdot \nabla )w_{ij} \psi _i (t),\;\overline{\overline {\omega }}_{i{m}'} )} } \\
 &\quad \;\quad \quad \quad \to -\int_0^\delta {\sum\limits_{i=1}^3 {((\overline
{u}\cdot \nabla )w_{ij} \psi _i (t),\;\overline{\overline {\omega
}}_i^\ast )} } = \int_0^\delta {\sum\limits_{i=1}^3 {((\overline
{u}\cdot \nabla ) \overline{\overline {\omega
}}_i^\ast ,\;w_{ij} \psi _i (t))} } \\
 &\int_0^\delta {\sum\limits_{i=1}^3 {((\overline{\overline {\omega }}_{i{m}'} \cdot \nabla
)\,\overline {u}_i ,\;w_{ij} \psi _i (t))} } \to \int_0^\delta
{\sum\limits_{i=1}^3 {((\overline{\overline {\omega} }^\ast \cdot
\nabla )\,\overline
{u}_i ,\;w_{ij} \psi _i (t))} }  \\
 &\sum\limits_{i=1}^3 {(\omega _{i0}^{{m}'} ,\;w_{ij} )\psi _i (0)} \to
\sum\limits_{i=1}^3 {(\omega _{i0},\;w_{ij} )\psi _i (0)}  \\
\end{split}
\end{equation*}
Thus, in the limit we find
\begin{equation}
\begin{split}
 &-\int_0^\delta {\sum\limits_{i=1}^3 {(\tilde\omega _i^\ast ,\;\partial _t \psi _i
(t)v_i )\,dt} } + \int_0^\delta {\sum\limits_{i=1}^3 {(\nabla \tilde\omega
_i^\ast
,\;\psi _i (t)\nabla v_i )\,dt} } \\
 &+ \int_0^\delta {\sum\limits_{i=1}^3 {((\overline {u}\cdot \nabla ) \overline{\overline {\omega
}}_i^\ast ,\;v_i \psi _i (t))dt} } - \int_0^\delta
{\sum\limits_{i=1}^3 {((\overline{\overline {\omega} }^\ast \cdot
\nabla )\,\overline {u}_i ,\;v_i \psi _i (t))} }
=\sum\limits_{i=1}^3 {(\omega _{i0} ,\;v_i )\psi _i (0)} \\
\end{split}
\end{equation}
holds for $v_i =w_{i1} ,\;w_{i2} ,\cdots $; by this equation holds
for $v_i =$ any finite linear combination of the $w_{ij} $, and by a
continuity argument above equation is still true for any $v_i \in
V$. Hence we find that $\tilde\omega _i^\ast (i=1,2,3)$ is a Leray-Hopf
weak solution of the system (9).
\\

Finally it remains to prove that $\tilde\omega _i^\ast $ satisfy the initial
conditions. For this we multiply (9) by $v_i \psi _i (t)$, after integrating
some terms by parts, we get in the same way,
\begin{equation*}
\begin{split}
 &-\int_0^\delta {\sum\limits_{i=1}^3 {(\tilde\omega _i^\ast ,\;\partial _t \psi _i
(t)v_i )} } + \int_0^\delta {\sum\limits_{i=1}^3 {(\nabla \tilde\omega _i^\ast
,\;\psi _i (t)\nabla v_i )\,dt} } \\
 &+ \int_0^\delta {\sum\limits_{i=1}^3 {((\overline {u}\cdot \nabla ) \overline{\overline {\omega
}}_i^\ast ,\;v_i \psi _i (t))} } - \int_0^\delta
{\sum\limits_{i=1}^3 {((\overline{\overline {\omega }}^\ast \cdot
\nabla )\,\overline {u}_i ,\;v_i \psi _i (t))} }
=\sum\limits_{i=1}^3 {(\tilde\omega _i^\ast (0),\;v_i )\psi _i (0)} \\
\end{split}
\end{equation*}
By comparison with (16),
\[
\sum\limits_{i=1}^3 {(\tilde\omega _i^\ast (0)-\omega _{i0} ,\;v_i )\psi _i (0)}
=0
\]
Therefore we can choose $\psi _i $ particularly such that
\[
(\tilde\omega _i^\ast (0)-\omega _{i0} ,\;v_i )=0,\quad \quad \forall\;
v_i \in V
\]
\\

\textbf{4. Convergence}

Now the partition is refined infinitely and $\varepsilon$ becomes
sufficiently small, we will prove that there exists some subsequence
of the solutions of auxiliary problems which converges to a weak
solution of (2).
\\

Since
\[
\mathop {\sup }\limits_{t\in (0,T)} \;\int_{{\mathbb R}^3} {(\tilde
{\omega }_1^2 +\tilde {\omega }_2^2 +\tilde {\omega }_3^2 )}
\;<+\infty
\]
the family $(\tilde {\omega }_1 ,\tilde {\omega }_2 ,\tilde {\omega
}_3 )$ is uniformly bounded in $L^2(0,T;H)\cap L^\infty (0,T;H)$,
then we can choose ${k}'\to \infty $ or $\Delta t_k ^\prime \to 0$
(in this case ${\varepsilon}' \to 0 $ and $ m' $ has to tend to $
\infty $), such that there exists a subsequence $({\tilde {\omega
}}'_1 ,{\tilde {\omega }}'_2 ,{\tilde {\omega }}'_3 )$ converging
weakly in $L^2(0,T;H)$ and weak-star in $L^\infty (0,T;H)$ to some
element $(\omega _1^\ast ,\omega _2^\ast ,\omega _3^\ast )$. On the
other hand, because $\tilde {\omega }_i (i=1,2,3)$ belong to
$L^2(0,T;H)$, we can verify that
\[
\overline {\omega }_i (x,t)=\left\{ {\frac{1}{\Delta t_k
}\int_{t_{k-1} }^{t_k } {\tilde {\omega }_i (x,t)dt} ,\;\;t\in
(t_{k-1} ,t_k )\subset (0,T)} \right\}
\]
also belongs to $L^2(0,T;H)$. In fact,
\begin{equation*}
\begin{split}
 &\int_0^T {\int_{{\mathbb R}^3} {\bar {\omega }_i^2 (x,t)} } =\sum\limits_k
{\int_{t_{k-1} }^{t_k } {\int_{{\mathbb R}^3} {\left(
{\frac{1}{\Delta t_k
}\int_{t_{k-1} }^{t_k } {\tilde {\omega }_i (x,t)} } \right)} } } ^2= \\
 &\quad \quad  =\sum\limits_k {\frac{1}{\Delta t_k^2 }\cdot
\Delta t_k \cdot \int_{{\mathbb R}^3} {\left( {\int_{t_{k-1} }^{t_k
} {\tilde {\omega }_i (x,t)} } \right)} } ^2\le \sum\limits_k
{\frac{1}{\Delta t_k }\int_{{\mathbb R}^3} {\int_{t_{k-1} }^{t_k } 1
\cdot \int_{t_{k-1} }^{t_k }
{\tilde {\omega }_i^2 (x,t)} } } \\
 &\quad \quad  =\sum\limits_k {\int_{t_{k-1} }^{t_k }
{\int_{{\mathbb R}^3} {\tilde {\omega }_i^2 (x,t)} } } =\int_0^T
{\int_{{\mathbb
R}^3} {\tilde {\omega }_i^2 (x,t)} } <+\infty \\
\end{split}
\end{equation*}

In the same way, we know from (6) that the function
\[
\overline {u}_i (x,t)=\left\{ {\frac{1}{\Delta t_k }\int_{t_{k-1}
}^{t_k } {u_i (x,t)dt} ,\;\;t\in (t_{k-1} ,t_k )\subset (0,T)}
\right\}
\]
belongs to $L^2(0,T;H)$.
\\

Finally we will prove that $(\omega _1^\ast ,\omega _2^\ast ,\omega _3^\ast
)$ is a solution of the vorticity-velocity form of Navier-Stokes equation
(2).

Taking $\varphi _i \in C^\infty ((0,T)\times {\mathbb
R}^3)\;\;(i=1,2,3)$, and
\[
\partial _{x_1 } \varphi _1 +\partial _{x_2 }
\varphi _2 +\partial _{x_3 } \varphi _3 =0
\]
we have
\begin{equation*}
\begin{split}
 &\sum\limits_{k=1}^N {\int_{t_{k-1} }^{t_k } {\int_{{\mathbb R}^3} {\theta
_{r} \varphi _1 (\partial _t \tilde {\omega }_1 \,+ \overline
{u}_1^k
\partial _{x_1 } \overline{\overline {\omega }}_1^k + \overline {u}_2^k \partial _{x_2 } \overline{\overline
{\omega }}_1^k + \overline {u}_3^k \partial _{x_3 } \overline{\overline {\omega }}_1^k -} } } \\
\end{split}
\end{equation*}
\begin{equation*}
\begin{split}
 &\quad \quad \quad \quad \quad \quad \quad \quad \quad - \overline{\overline {\omega }}_1^k
\partial _{x_1 } \overline {u}_1^k - \overline{\overline {\omega }}_2^k \partial _{x_2 } \overline
{u}_1^k - \overline{\overline {\omega }}_3^k \partial _{x_3 }
\overline {u}_1^k + \partial _{x_1 }
q - \Delta \tilde {\omega }_1 ) = 0 \\
 &\sum\limits_{k=1}^N {\int_{t_{k-1} }^{t_k } {\int_{{\mathbb R}^3} {\theta
_{r} \varphi _2 (\partial _t \tilde {\omega }_2 + \overline {u}_1^k
\partial _{x_1 } \overline{\overline {\omega }}_2^k + \overline {u}_2^k \partial _{x_2 } \overline{\overline
{\omega }}_2^k + \overline {u}_3^k \partial _{x_3 } \overline{\overline {\omega }}_2^k } } } - \\
 &\quad \quad \quad \quad \quad \quad \quad \quad \quad - \overline{\overline {\omega }}_1^k
\partial _{x_1 } \overline {u}_2^k - \overline{\overline {\omega }}_2^k \partial _{x_2 } \overline
{u}_2^k - \overline{\overline {\omega }}_3^k \partial _{x_3 }
\overline {u}_2^k + \partial _{x_2 }
q - \Delta \tilde {\omega }_2 ) = 0 \\
 &\sum\limits_{k=1}^N {\int_{t_{k-1} }^{t_k } {\int_{{\mathbb R}^3} {\theta
_{r} \varphi _3 (\partial _t \tilde {\omega }_3 + \overline {u}_1^k
\partial _{x_1 } \overline{\overline {\omega }}_3^k + \overline {u}_2^k \partial _{x_2 } \overline{\overline
{\omega }}_3^k + \overline {u}_3^k \partial _{x_3 } \overline{\overline {\omega} }_3^k } } } - \\
 &\quad \quad \quad \quad \quad \quad \quad \quad \quad - \overline{\overline {\omega }}_1^k
\partial _{x_1 } \overline {u}_3^k - \overline{\overline {\omega }}_2^k \partial _{x_2 } \overline
{u}_3^k - \overline{\overline {\omega }}_3^k \partial _{x_3 }
\overline {u}_3^k + \partial _{x_3 }
q - \Delta \tilde {\omega }_3 ) = 0 \\
\end{split}
\end{equation*}
Here $\tilde {\omega }_i \;(i=1,2,3)$ denote the collection of those
solutions of problem (7) defined on every $(t_{k-1},t_k)$.
Integrating by parts we get
\begin{equation*}
\begin{split}
 &\sum\limits_{k=1}^N {\int_{t_{k-1} }^{t_k } {\int_{{\mathbb R}^3} {\theta
_{r} (\tilde {\omega }_1 \partial _t \varphi _1 \,+
\overline{\overline {\omega} }_1^k ((\overline {u}_1^k \partial
_{x_1 } \varphi _1 +\varphi _1 \,\partial _{x_1 } \overline {u}_1^k
) + (\overline {u}_2^k
\partial _{x_2 } \varphi _1 +\varphi _1
\,\partial _{x_2 } \overline {u}_2^k )+} } } \\
 &\quad + (\overline {u}_3^k \partial _{x_3 } \varphi _1 +\varphi _1 \,\partial
_{x_3 } \overline {u}_3^k )) - \overline {u}_1^k
((\overline{\overline {\omega} }_1^k
\partial _{x_1 } \varphi _1 +\varphi _1 \,\partial _{x_1 } \overline{\overline
{\omega }}_1^k ) + (\overline{\overline {\omega }}_2^k \partial
_{x_2 } \varphi _1 +\varphi _1 \,\partial _{x_2 }
\overline{\overline
{\omega }}_2^k ) + \\
 &\quad + (\overline{\overline {\omega }}_3^k \partial _{x_3 } \varphi _1 +\varphi _1
\,\partial _{x_3 } \overline{\overline {\omega }}_3^k )) + q
\partial _{x_1 } \varphi _1
+\tilde {\omega }_1 \Delta \varphi _1 ) + \\
 &+ \sum\limits_{k=1}^N {\int_{t_{k-1} }^{t_k } {\int_{{\mathbb R}^3} {(\overline{\overline
{\omega }}_1^k (\varphi _1 \overline {u}_1^k \partial _{x_1 } \theta
_{r} + \varphi _1 \overline {u}_2^k \partial _{x_2 } \theta _{r}
+\varphi _1
\overline {u}_3^k \partial _{x_3 } \theta _{r} )-} } } \\
 &\quad - \overline {u}_1^k (\varphi _1 \overline{\overline {\omega }}_1^k \partial _{x_1 } \theta
_{r} + \varphi _1 \overline{\overline {\omega }}_2^k \partial _{x_2
} \theta _{r} + \varphi _1 \overline{\overline {\omega }}_3^k
\partial _{x_3 } \theta
_{r} ) \\
 &\quad + q \varphi _1 \partial _{x_1 } \theta _{r} +\tilde {\omega
}_1 \varphi _1 \Delta \theta _{r} +2\tilde {\omega }_1 (\partial
_{x_1 } \theta _{r} \partial _{x_1 } \varphi _1 +\partial _{x_2 }
\theta _{r} \partial _{x_2 } \varphi _1 +\partial _{x_3 } \theta
_{r} \partial _{x_3 } \varphi _1 )) \\
 &\quad =\sum\limits_{k=1}^N {\int_{{\mathbb R}^3} {\theta _{r} (\varphi
_1 (x,t_k )\tilde {\omega }_1 (x,t_k ) - \varphi _1 (x,t_{k-1}
)\tilde {\omega
}_1 (x,t_{k-1} ))} } \\
 &\sum\limits_{k=1}^N {\int_{t_{k-1} }^{t_k } {\int_{{\mathbb R}^3} {\theta
_{r} (\tilde {\omega }_2 \partial _t \varphi _2 \,+
\overline{\overline {\omega }}_2^k ((\overline {u}_1^k \partial
_{x_1 } \varphi _2 +\varphi _2 \,\partial _{x_1 } \overline {u}_1^k
) + (\overline {u}_2^k
\partial _{x_2 } \varphi _2 +\varphi _2
\,\partial _{x_2 } \overline {u}_2^k )+} } } \\
 &\quad + (\overline {u}_3^k \partial _{x_3 } \varphi _2 +\varphi _2 \,\partial
_{x_3 } \overline {u}_3^k )) - \overline {u}_2^k
((\overline{\overline {\omega }}_1^k
\partial _{x_1 } \varphi _2 +\varphi _2 \,\partial _{x_1 } \overline{\overline
{\omega }}_1^k ) + (\overline{\overline {\omega }}_2^k \partial
_{x_2 } \varphi _2 +\varphi _2 \,\partial _{x_2 }
\overline{\overline
{\omega }}_2^k )+ \\
 &\quad + (\overline{\overline {\omega }}_3^k \partial _{x_3 } \varphi _2 +\varphi _2
\,\partial _{x_3 } \overline{\overline {\omega }}_3^k )) + q
\partial _{x_2 } \varphi _2
+\tilde {\omega }_2 \Delta \varphi _2 )+ \\
 &+\sum\limits_{k=1}^N {\int_{t_{k-1} }^{t_k } {\int_{{\mathbb R}^3} {(\overline{\overline
{\omega }}_2^k (\varphi _2 \overline {u}_1^k \partial _{x_1 } \theta
_{r} + \varphi _2 \overline {u}_2^k \partial _{x_2 } \theta _{r}
+\varphi _2
\overline {u}_3^k \partial _{x_3 } \theta _{r} )-} } } \\
 &\quad - \overline {u}_2^k (\varphi _2 \overline{\overline {\omega }}_1^k \partial _{x_1 } \theta
_{r} + \varphi _2 \overline{\overline {\omega }}_2^k \partial _{x_2
} \theta _{r} + \varphi _2 \overline{\overline {\omega }}_3^k
\partial _{x_3 } \theta
_{r} ) \\
 &\quad + q\varphi _2 \partial _{x_2 } \theta _{r} + \tilde {\omega
}_2 \varphi _2 \Delta \theta _{r} + 2\tilde {\omega }_2 (\partial
_{x_1 } \theta _{r} \partial _{x_1 } \varphi _2 + \partial _{x_2 }
\theta _{r} \partial _{x_2 } \varphi _2 + \partial _{x_3 } \theta
_{r} \partial _{x_3 } \varphi _2 )) \\
 &\quad =\sum\limits_{k=1}^N {\int_{{\mathbb R}^3} {\theta _{r} (\varphi
_2 (x,t_k )\tilde {\omega }_2 (x,t_k ) - \varphi _2 (x,t_{k-1}
)\tilde {\omega
}_2 (x,t_{k-1} ))} } \\
 &\sum\limits_{k=1}^N {\int_{t_{k-1} }^{t_k } {\int_{{\mathbb R}^3} {\theta
_{r} (\tilde {\omega }_3 \partial _t \varphi _3 \,+
\overline{\overline {\omega }}_3^k ((\overline {u}_1^k \partial
_{x_1 } \varphi _3 +\varphi _3 \,\partial _{x_1 } \overline {u}_1^k
) + (\overline {u}_2^k
\partial _{x_2 } \varphi _3 + \varphi _3
\,\partial _{x_2 } \overline {u}_2^k )+ } } } \\
 &\quad + (\overline {u}_3^k \partial _{x_3 } \varphi _3 + \varphi _3 \,\partial
_{x_3 } \overline {u}_3^k )) - \overline {u}_3^k
((\overline{\overline {\omega }}_1^k
\partial _{x_1 } \varphi _3 +\varphi _3 \,\partial _{x_1 } \overline{\overline
{\omega }}_1^k ) + (\overline{\overline {\omega} }_2^k \partial
_{x_2 } \varphi _3 +\varphi _3 \,\partial _{x_2 }
\overline{\overline
{\omega }}_2^k ) + \\
\end{split}
\end{equation*}
\begin{equation*}
\begin{split}
 &\quad + (\overline{\overline {\omega }}_3^k \partial _{x_3 } \varphi _3 + \varphi _3
\,\partial _{x_3 } \overline{\overline {\omega }}_3^k )) + q
\partial _{x_3 } \varphi _3
+ \tilde {\omega }_3 \Delta \varphi _3 ) + \\
 &+ \sum\limits_{k=1}^N {\int_{t_{k-1} }^{t_k } {\int_{{\mathbb R}^3} {(\overline{\overline
{\omega }}_3^k (\varphi _3 \overline {u}_1^k \partial _{x_1 } \theta
_{r} + \varphi _3 \overline {u}_2^k \partial _{x_2 } \theta _{r} +
\varphi _3
\overline {u}_3^k \partial _{x_3 } \theta _{r} )-} } } \\
 &\quad - \overline {u}_3^k (\varphi _3 \overline{\overline {\omega }}_1^k \partial _{x_1 } \theta
_{r} + \varphi _3 \overline{\overline {\omega} }_2^k \partial _{x_2
} \theta _{r} + \varphi _3 \overline{\overline {\omega }}_3^k
\partial _{x_3 } \theta
_{r} ) \\
 &\quad + q \varphi _3 \partial _{x_3 } \theta _{r} + \tilde {\omega
}_3 \varphi _3 \Delta \theta _{r} + 2 \tilde {\omega }_3 (\partial
_{x_1 } \theta _{r} \partial _{x_1 } \varphi _3 + \partial _{x_2 }
\theta _{r} \partial _{x_2 } \varphi _3 + \partial _{x_3 } \theta
_{r} \partial _{x_3 } \varphi _3 )) \\
 &\quad = \sum\limits_{k=1}^N {\int_{{\mathbb R}^3} {\theta _{r} (\varphi
_3 (x,t_k )\tilde {\omega }_3 (x,t_k ) - \varphi _3 (x,t_{k-1}
)\tilde {\omega
}_3 (x,t_{k-1} ))} } \\
\end{split}
\end{equation*}
Let $ r\to +\infty $,
\begin{equation*}
\begin{split}
 &\sum\limits_{k=1}^N {\int_{t_{k-1} }^{t_k } {\int_{{\mathbb R}^3} {(\tilde
{\omega }_1 \partial _t \varphi _1 \,+ \overline{\overline {\omega
}}_1^k ((\overline {u}_1^k
\partial _{x_1 } \varphi _1 +\varphi _1 \,\partial _{x_1 } \overline {u}_1^k
) + (\overline {u}_2^k \partial _{x_2 } \varphi _1 + \varphi _1
\,\partial _{x_2 }
\overline {u}_2^k )+} } } \\
 &\quad +(\overline {u}_3^k \partial _{x_3 } \varphi _1 + \varphi _1 \,\partial
_{x_3 } \overline {u}_3^k )) - \overline {u}_1^k
((\overline{\overline {\omega} }_1^k
\partial _{x_1 } \varphi _1 + \varphi _1 \,\partial _{x_1 } \overline{\overline
{\omega} }_1^k ) + (\overline{\overline {\omega }}_2^k \partial
_{x_2 } \varphi _1 + \varphi _1 \,\partial _{x_2 }
\overline{\overline
{\omega} }_2^k )+ \\
 &\quad + (\overline{\overline {\omega }}_3^k \partial _{x_3 } \varphi _1 +\varphi _1
\,\partial _{x_3 } \overline{\overline {\omega }}_3^k )) + q
\partial _{x_1 } \varphi _1
+ \tilde {\omega }_1 \Delta \varphi _1 ) \\
 &\quad = \sum\limits_{k=1}^N {\int_{{\mathbb R}^3} {(\varphi _1 (x,t_k )\tilde
{\omega }_1 (x,t_k )-\varphi _1 (x,t_{k-1} )\tilde {\omega }_1
(x,t_{k-1}
))} } \\
 &\sum\limits_{k=1}^N {\int_{t_{k-1} }^{t_k } {\int_{{\mathbb R}^3} {(\tilde
{\omega }_2 \partial _t \varphi _2 \,+ \overline{\overline {\omega
}}_2^k ((\overline {u}_1^k
\partial _{x_1 } \varphi _2 + \varphi _2 \,\partial _{x_1 } \overline {u}_1^k
) + (\overline {u}_2^k \partial _{x_2 } \varphi _2 + \varphi _2
\,\partial _{x_2 }
\overline {u}_2^k )+ } } } \\
 &\quad + (\overline {u}_3^k \partial _{x_3 } \varphi _2 +\varphi _2 \,\partial
_{x_3 } \overline {u}_3^k )) - \overline {u}_2^k
((\overline{\overline {\omega }}_1^k
\partial _{x_1 } \varphi _2 + \varphi _2 \,\partial _{x_1 } \overline{\overline
{\omega} }_1^k ) + (\overline{\overline {\omega }}_2^k \partial
_{x_2 } \varphi _2 + \varphi _2 \,\partial _{x_2 }
\overline{\overline
{\omega }}_2^k ) + \\
 &\quad + (\overline{\overline {\omega} }_3^k \partial _{x_3 } \varphi _2 + \varphi _2
\,\partial _{x_3 } \overline{\overline {\omega} }_3^k )) + q
\partial _{x_2 } \varphi _2
+ \tilde {\omega }_2 \Delta \varphi _2 ) \\
 &\quad = \sum\limits_{k=1}^N {\int_{{\mathbb R}^3} {(\varphi _2 (x,t_k )\tilde
{\omega }_2 (x,t_k )-\varphi _2 (x,t_{k-1} )\tilde {\omega }_2 (x,t_{k-1}
))} } \\
 &\sum\limits_{k=1}^N {\int_{t_{k-1} }^{t_k } {\int_{{\mathbb R}^3} {(\tilde
{\omega }_3 \partial _t \varphi _3 \,+ \overline{\overline {\omega
}}_3^k ((\overline {u}_1^k
\partial _{x_1 } \varphi _3 + \varphi _3 \,\partial _{x_1 } \overline {u}_1^k
) + (\overline {u}_2^k \partial _{x_2 } \varphi _3 + \varphi _3
\,\partial _{x_2 }
\overline {u}_2^k )+} } } \\
 &\quad + (\overline {u}_3^k \partial _{x_3 } \varphi _3 + \varphi _3 \,\partial
_{x_3 } \overline {u}_3^k )) - \overline {u}_3^k
((\overline{\overline {\omega }}_1^k
\partial _{x_1 } \varphi _3 + \varphi _3 \,\partial _{x_1 } \overline{\overline
{\omega} }_1^k ) + (\overline{\overline {\omega }}_2^k \partial
_{x_2 } \varphi _3 + \varphi _3 \,\partial _{x_2 }
\overline{\overline
{\omega} }_2^k )+ \\
 &\quad + (\overline{\overline {\omega }}_3^k \partial _{x_3 } \varphi _3 + \varphi _3
\,\partial _{x_3 } \overline{\overline {\omega }}_3^k )) + q
\partial _{x_3 } \varphi _3
+ \tilde {\omega }_3 \Delta \varphi _3 ) \\
 &\quad = \sum\limits_{k=1}^N {\int_{{\mathbb R}^3} {(\varphi _3 (x,t_k )\tilde
{\omega }_3 (x,t_k )-\varphi _3 (x,t_{k-1} )\tilde {\omega }_3 (x,t_{k-1}
))} } \\
\end{split}
\end{equation*}
\\

For a certain solution $u$ of (1), we can prove due to (6) that
\\

$\overline {u}_i \to u_i $ \textit{in} $L^2(0,T;H)$ \textit{strongly}
\\

\noindent as $k \to \infty $, or $\Delta t_k \to 0$.
\\

In fact, set $Q=(0,T)\times {\mathbb R}^3$, $\Delta t=\mathop {\max
}\limits_k \{\Delta t_k \}$, $\forall \varepsilon >0$, and $u_i \in
L^2(0,T;L^2({\mathbb R}^3))$, there exists a $v_i \in C^\infty
(0,T;L^2({\mathbb R}^3))$ such that
\[
\left\| {\,u_i -v_i } \right\|_{L^2(Q)} <\varepsilon
\]
By means of the same partition as that for $\overline {u}_i $ to
construct $\overline {v}_i $, since there exists a constant $C>0$
such that $\left\| {\partial _t v_i } \right\|_{L^2({\mathbb R}^3)}
\le C$, and $\mathop {\max }\limits_t \left\| {\,\overline {v}_i
-v_i } \right\|_{L^2({\mathbb R}^3)} \le C\;\Delta t$, it follows
that
\[
\left\| {\,\overline {v}_i -v_i } \right\|_{L^2(Q)} =\left(
{\int_0^T {\left\| {\,\overline {v}_i -v_i } \right\|_{L^2({\mathbb
R}^3)}^2 } } \right)^{1/2}\le C\,T^{1/2}\Delta t
\]
Thus
\[
\overline {v}_i \to v_i \quad \left( {\;L^\infty (0,T;L^2({\mathbb
R}^3))\;} \right),\quad \mbox{as}\;\,\Delta t\to 0
\]
Take $\Delta t$ such that $\left\| {\overline {v}_i -v_i }
\right\|_{L^2(Q)} <\varepsilon $. Moreover,
\begin{equation*}
\begin{split}
 &\int_0^T {\left\| {\,\overline {u}_i - \overline {v}_i } \right\|_{L^2({\mathbb R}^3)}^2 }
=\sum\limits_{k=1}^N {\left\| {\frac{1}{\Delta t_k }\int_{t_{k-1}
}^{t_k }
{(u_i -v_i )} } \right\|} _{L^2({\mathbb R}^3)}^2 \Delta t_k \\
 &\quad \le \sum\limits_{k=1}^N {\left\| {\;\left( {\int_{t_{k-1} }^{t_k }
{(u_i -v_i )^2} } \right)^{1/2}} \right\|} _{L^2({\mathbb R}^3)}^2
\le \int_0^T
{\left\| {\,u_i -v_i } \right\|_{L^2({\mathbb R}^3)}^2 } \\
\end{split}
\end{equation*}
so that $\left\| {\,\overline {u}_i - \overline {v}_i }
\right\|_{L^2(Q)} \le \left\| {u_i -v_i } \right\|_{L^2(Q)}
<\varepsilon $. Therefore,
\[
\left\| {\,\overline {u}_i -u_i } \right\|_{L^2(Q)} \le \left\|
{\,u_i -v_i } \right\|_{L^2(Q)} + \left\| {\,v_i - \overline {v}_i }
\right\|_{L^2(Q)} +\left\| {\,\overline {v}_i - \overline {u}_i }
\right\|_{L^2(Q)} <3\varepsilon
\]
Hence as $\Delta t\to 0$, we have $\left\| {\,\overline {u}_i -u_i }
\right\|_{L^2(Q)} \to 0$.
\\

On the other hand, from section 2 we have the following conclusions :
\\

$\tilde {\omega }_i^\prime \to \omega _i^\ast $ \textit{in} $L^2
(0,T;H)$ \textit{weakly}, \textit{and in} $L^\infty (0,T;H)$
\textit{weak-star}
\\

\noindent for a subsequence as ${k^\prime} \to \infty $, or $\Delta
t_k ^\prime \to 0$.
\\

For a $ {\omega }_i \in L^2 (0,T;H) $, similar to above we know that
\\

$\overline {\omega }_i \to \omega _i $ \textit{in} $L^2(0,T;H)$
\textit{strongly}
\\

\noindent as ${k} \to \infty $ or $\Delta t_k \to 0$.
\\

Moreover, set $B_\varepsilon = \{x:\left| x \right| < \varepsilon \}$, then
\begin{equation*}
\begin{split}
 &\left\| \,{\overline{\overline \omega} _i - \overline \omega_i } \right\|_{L^2({{\mathbb R}^3} )} =
\left\| \;\; {\int\limits_{\left| y \right| \le \varepsilon }
{J_\varepsilon (y)\, [\,\overline \omega _i (x - y) - \overline
\omega _i (x)\,]\, dy} } \right\|_{L^2({{\mathbb R}^3}
)} \\
\end{split}
\end{equation*}
\begin{equation*}
\begin{split}
 &\qquad \le \int\limits_{\left| y \right| \le \varepsilon } {J_\varepsilon
(y)\left\| \, {\overline \omega _i (x - y) - \overline \omega _i
(x)}
\right\|_{L^2({{\mathbb R}^3} )} dy} \\
 &\qquad \le  \mathop {\sup }\limits_{\left| y
\right| \le \varepsilon } \left\| \, {\overline \omega _i (x - y) -
\overline \omega _i (x)}
\right\|_{L^2({{\mathbb R}^3} )} \to 0 \\
 \end{split}
\end{equation*}
as $\varepsilon \to 0$, and
\[
\left\| \, {\overline{\overline \omega} _i - \omega_i }
\right\|_{L^2(Q)} \le \left( {\int_0^T {\left\| \,
{\overline{\overline \omega} _i - \overline \omega _i }
\right\|_{L^2({{\mathbb R}^3} )}^2 } } \right)^{1 / 2} +\; \left\|
\, {\overline \omega _i - \omega_i } \right\|_{L^2(Q)} \to 0
\]
as ${k} \to \infty $ or $\Delta t_k \to 0$.
\\

Thus we finally we obtain
\[
{\begin{array}{*{20}c}
 {\overline {\omega }_i^{k'} \to {\omega} _i^{*} } \hfill \\
 {\overline {\overline {\omega }}_i^{k'} \to \omega _i^{*} } \hfill \\
\end{array} }\quad \quad \mbox{\textit{in}}\;\;L^2(0,T,H)\;\;\mbox{\textit{weakly}}
\]
as ${k^\prime} \to \infty $, or $\Delta t_k ^\prime \to 0$.
\\

These convergence results enable us to pass the limit. That is,
\begin{equation*}
\begin{split}
 &\sum\limits_{{k}'} {\int_{t_{{k}'-1} }^{t_{{k}'} } {\int_{{\mathbb R}^3}
{(\tilde {\omega }_1 \partial _t \varphi _1 \,+ \overline{\overline
{\omega }}_1^{{k}'} (\overline {u}_1^{{k}'} \partial _{x_1 } \varphi
_1 + \overline {u}_2^{{k}'}
\partial _{x_2 }
\varphi _1 + \overline {u}_3^{{k}'} \partial _{x_3 } \varphi _1 )-} } } \\
 &\quad \quad \quad \quad \quad \quad - \overline {u}_1^{{k}'} (\overline{\overline {\omega
}}_1^{{k}'} \partial _{x_1 } \varphi _1 + \overline{\overline
{\omega }}_2^{{k}'}
\partial _{x_2 } \varphi _1 + \overline{\overline {\omega }}_3^{{k}'} \partial _{x_3
} \varphi _1
) + q \partial _{x_1 } \varphi _1 + \tilde {\omega }_1 \Delta \varphi _1 ) \\
 &\quad \quad \quad \quad \quad \quad = \int_{{\mathbb R}^3} {(\varphi _1
(x,T)\tilde {\omega }_1 (x,T)-\varphi _1 (x,0) \tilde {\omega }_1 (x,0))} \\
 &\sum\limits_{{k}'} {\int_{t_{{k}'-1} }^{t_{{k}'} } {\int_{{\mathbb R}^3}
{(\tilde {\omega }_2 \partial _t \varphi _2 \,+ \overline{\overline
{\omega }}_2^{{k}'} (\overline {u}_1^{{k}'} \partial _{x_1 } \varphi
_2 + \overline {u}_2^{{k}'} \partial _{x_2 }
\varphi _2 + \overline {u}_3^{{k}'} \partial _{x_3 } \varphi _2 )-} } } \\
 &\quad \quad \quad \quad \quad \quad - \overline {u}_2^{{k}'} (\overline{\overline {\omega
}}_1^{{k}'} \partial _{x_1 } \varphi _2 + \overline{\overline
{\omega }}_2^{{k}'}
\partial _{x_2 } \varphi _2 + \overline{\overline {\omega }}_3^{{k}'} \partial _{x_3
} \varphi _2
) + q \partial _{x_2 } \varphi _2 + \tilde {\omega }_2 \Delta \varphi _2 ) \\
 &\quad \quad \quad \quad \quad \quad = \int_{{\mathbb R}^3} {(\varphi _2
(x,T)\tilde {\omega }_2 (x,T) - \varphi _2 (x,0)\tilde {\omega }_2 (x,0))} \\
 &\sum\limits_{{k}'} {\int_{t_{{k}'-1} }^{t_{{k}'} } {\int_{{\mathbb R}^3}
{(\tilde {\omega }_3 \partial _t \varphi _3 \,+ \overline{\overline
{\omega }}_3^{{k}'} (\overline {u}_1^{{k}'} \partial _{x_1 } \varphi
_3 + \overline {u}_2^{{k}'} \partial _{x_2 }
\varphi _3 + \overline {u}_3^{{k}'} \partial _{x_3 } \varphi _3 )-} } } \\
 &\quad \quad \quad \quad \quad \quad - \overline {u}_3^{{k}'} (\overline{\overline {\omega
}}_1^{{k}'} \partial _{x_1 } \varphi _3 + \overline{\overline
{\omega} }_2^{{k}'}
\partial _{x_2 } \varphi _3 + \overline{\overline {\omega }}_3^{{k}'} \partial _{x_3
} \varphi _3
) + q \partial _{x_3 } \varphi _3 + \tilde {\omega }_3 \Delta \varphi _3 ) \\
 &\quad \quad \quad \quad \quad \quad = \int_{{\mathbb R}^3} {(\varphi _3
(x,T)\tilde {\omega }_3 (x,T) - \varphi _3 (x,0)\tilde {\omega }_3 (x,0))} \\
\end{split}
\end{equation*}
\\
This is equivalent to
\begin{equation*}
\begin{split}
 &\int_0^T {\int_{{\mathbb R}^3} {\,\{(\omega _1^\ast \partial _t \varphi _1
\,+ \omega _2^\ast \partial _t \varphi _2 \,+ \omega _3^\ast
\partial _t
\varphi _3 )+ } } \\
 &\qquad \quad  + (\omega _1^\ast \Delta \varphi _1 + \omega _2^\ast \Delta
\varphi _2 + \omega _3^\ast \Delta \varphi _3 )+ \\
\end{split}
\end{equation*}
\begin{equation*}
\begin{split}
 &\quad + \omega _1^\ast (u_1 \partial _{x_1 } \varphi _1 + u_2 \partial _{x_2
} \varphi _1 + u_3 \partial _{x_3 } \varphi _1 ) + \omega _2^\ast
(u_1
\partial _{x_1 } \varphi _2 + u_2 \partial _{x_2 } \varphi _2 + u_3
\partial _{x_3 }
\varphi _2 ) + \\
 &\quad + \omega _3^\ast (u_1 \partial _{x_1 } \varphi _3 + u_2 \partial _{x_2
} \varphi _3 + u_3 \partial _{x_3 } \varphi _3 ) \\
 &\quad - u_1 (\omega _1^\ast \partial _{x_1 } \varphi _1 + \omega _2^\ast
\partial _{x_2 } \varphi _1 + \omega _3^\ast \partial _{x_3 } \varphi _1
) - u_2 (\omega _1^\ast \partial _{x_1 } \varphi _2 + \omega _2^\ast
\partial
_{x_2 } \varphi _2 + \omega _3^\ast \partial _{x_3 } \varphi _2 ) - \\
 &\quad - u_3 (\omega _1^\ast \partial _{x_1 } \varphi _3 + \omega _2^\ast
\partial _{x_2 } \varphi _3 + \omega _3^\ast \partial _{x_3 } \varphi _3 )\}
\\
 &=\int_{{\mathbb R}^3} {\{(\varphi _1 (x,T)\omega _1^\ast (x,T) + \varphi _2
(x,T)\omega _2^\ast (x,T) + \varphi _3 (x,T)\omega _3^\ast (x,T))-} \\
 &\quad \quad \;\;- (\varphi _{10} (x)\omega _{10} (x) + \varphi _{20} (x)\omega
_{20} (x) + \varphi _{30} (x)\omega _{30} (x))\} \\
\end{split}
\end{equation*}
\\
Here we also have
\[
\omega _i^\ast (x,0)=\omega _{i0} (x),\quad \varphi _i (x,0)=\varphi _{i0}
(x),\quad i=1,2,3
\]
Hence we know that there exists some $\omega _i^\ast $ which belongs
to $L^\infty (0,T;L^2({\mathbb R}^3))$ and is a Leray-Hopf weak
solution of (2).
\\

Note that a weak formulation of the following equations:
\begin{equation*}
\begin{split}
 &\;\; \omega = \mbox{curl}\,u   \\
 &\int_0^{T} {\int_{{\mathbb R}^3} {\;\varphi \cdot [\partial _t \omega + (u \cdot
\nabla )\omega - (\omega \cdot \nabla )u - \nu \Delta \omega } } ] =
0
\\
\end{split}
\end{equation*}
is equivalent to
\[
\int_0^{T} {\int_{{\mathbb R}^3} {\;\tilde {\varphi } \cdot
[\partial _t u + (u \cdot \nabla )u + \nabla p - \nu \Delta u} } ] =
0
\]
for any $\varphi \in C^\infty ((0,T)\times {\mathbb{R}}^3)$ with
$\varphi_i \in C_0^\infty(\Omega)$ and zero extension outside
$\Omega, \forall\,\Omega \subset {\mathbb R}^3$, and $\tilde
{\varphi } = \mbox{curl}\varphi $, in some distribution sense.
\\
\\
\\
\\

\textbf{5. Regularity}

We can still use Galerkin procedure as in Section 2 and 3. Since $V$ is separable
there exists a sequence of linearly independent elements $w_{i1} ,\;\cdots
,\;w_{im} ,\;\cdots $ which is total in $V$. For each $m$ we define an
approximate solution $u_{im} $ of (1) as follows:
\[
u_{im} =\sum\limits_{j=1}^m {g_{ij} (t)\;w_{ij} }
\]
and by means of weighted function $\theta _{r} $
\begin{equation*}
\begin{split}
 &\int_{{\mathbb R}^3} {\theta _{r} w_{1j} \partial _t u_{1m} }
+\int_{{\mathbb R}^3} {\theta _{r} (u_{1m} \partial _{x_1 } u_{1m}
+u_{2m} \partial _{x_2 } u_{1m} +u_{3m} \partial _{x_3 } u_{1m}
)w_{1j} } +
\\
 &\quad \quad \quad \quad \quad \quad \quad +\int_{{\mathbb R}^3} {\theta
_{r} w_{1j} \partial _{x_1 } p} =\int_{{\mathbb R}^3}
{\theta _{r} w_{1j} \Delta u_{1m} } \\
\end{split}
\end{equation*}
\begin{equation}
\begin{split}
 &\int_{{\mathbb R}^3} {\theta _{r} w_{2j} \partial _t u_{2m} }
+\int_{{\mathbb R}^3} {\theta _{r} (u_{1m} \partial _{x_1 } u_{2m}
+u_{2m} \partial _{x_2 } u_{2m} +u_{3m} \partial _{x_3 } u_{2m}
)w_{2j} } +
\\
 &\quad \quad \quad \quad \quad \quad \quad +\int_{{\mathbb R}^3} {\theta
_{r} w_{2j} \partial _{x_2 } p} =\int_{{\mathbb R}^3} {\theta
_{r} w_{2j} \Delta u_{2m} } \\
 &\int_{{\mathbb R}^3} {\theta _{r} w_{3j} \partial _t u_{3m} }
+\int_{{\mathbb R}^3} {\theta _{r} (u_{1m} \partial _{x_1 } u_{3m}
+u_{2m} \partial _{x_2 } u_{3m} +u_{3m} \partial _{x_3 } u_{3m}
)w_{3j} } +
\\
 &\quad \quad \quad \quad \quad \quad \quad +\int_{{\mathbb R}^3} {\theta
_{r} w_{3j} \partial _{x_3 } p} =\int_{{\mathbb R}^3} {\theta
_{r} w_{3j} \Delta u_{3m} } \\
 &\quad \quad u_{im} (0)=u_{i0}^m ,\quad \quad j=1,\cdots ,m \\
\end{split}
\end{equation}
where $u_{i0}^m $ is the orthogonal projection in $H$ of $u_{i0} $
on the space spanned by $w_{i1} ,\;\cdots ,\;w_{im} $.
\\

We now are allowed to differentiate (17) in the $t$, we get
\begin{equation}
\begin{split}
 &\int_{{\mathbb R}^3} {\theta _{r} w_{1j} \partial _t^2 u_{1m} }
+\int_{{\mathbb R}^3} {\theta _{r} (\partial _t u_{1m}
\partial _{x_1 } u_{1m} +\partial _t u_{2m} \partial _{x_2 } u_{1m}
+\partial _t u_{3m}
\partial _{x_3 } u_{1m} )w_{1j} } + \\
 &\quad \quad \quad \quad \quad \quad \quad +\int_{{\mathbb R}^3} {\theta
_{r} (u_{1m} \partial _{x_1 } \partial _t u_{1m} +u_{2m}
\partial _{x_2 } \partial _t u_{1m} +u_{3m} \partial _{x_3 }
\partial _t u_{1m}
)w_{1j} } + \\
 &\quad \quad \quad \quad \quad \quad \quad +\int_{{\mathbb R}^3} {\theta
_{r} w_{1j} \partial _{x_1 } \partial _t p} =\int_{{\mathbb R}^3}
{\theta _{r} w_{1j} \Delta \partial _t u_{1m} } \\
 &\int_{{\mathbb R}^3} {\theta _{r} w_{2j} \partial _t^2 u_{2m} }
+\int_{{\mathbb R}^3} {\theta _{r} (\partial _t u_{1m}
\partial _{x_1 } u_{2m} +\partial _t u_{2m} \partial _{x_2 } u_{2m}
+\partial _t u_{3m}
\partial _{x_3 } u_{2m} )w_{2j} } + \\
 &\quad \quad \quad \quad \quad \quad \quad +\int_{{\mathbb R}^3} {\theta
_{r} (u_{1m} \partial _{x_1 } \partial _t u_{2m} +u_{2m}
\partial _{x_2 } \partial _t u_{2m} +u_{3m} \partial _{x_3 }
\partial _t u_{2m}
)w_{2j} } + \\
 &\quad \quad \quad \quad \quad \quad \quad +\int_{{\mathbb R}^3} {\theta
_{r} w_{2j} \partial _{x_2 } \partial _t p} =\int_{{\mathbb R}^3}
{\theta _{r} w_{2j} \Delta \partial _t u_{2m} } \\
 &\int_{{\mathbb R}^3} {\theta _{r} w_{3j} \partial _t^2 u_{3m} }
+\int_{{\mathbb R}^3} {\theta _{r} (\partial _t u_{1m}
\partial _{x_1 } u_{3m} +\partial _t u_{2m} \partial _{x_2 } u_{3m}
+\partial _t u_{3m}
\partial _{x_3 } u_{3m} )w_{3j} } + \\
 &\quad \quad \quad \quad \quad \quad \quad +\int_{{\mathbb R}^3} {\theta
_{r} (u_{1m} \partial _{x_1 } \partial _t u_{3m} +u_{2m}
\partial _{x_2 } \partial _t u_{3m} +u_{3m} \partial _{x_3 }
\partial _t u_{3m}
)w_{3j} } + \\
 &\quad \quad \quad \quad \quad \quad \quad +\int_{{\mathbb R}^3} {\theta
_{r} w_{3j} \partial _{x_3 } \partial _t p} =\int_{{\mathbb R}^3}
{\theta _{r} w_{3j} \Delta \partial _t u_{3m} } \\
 &\quad \quad \quad \quad \quad \quad \quad \quad \quad \quad \quad \quad
\quad j=1,\cdots ,m \\
\end{split}
\end{equation}
We multiply (18) by ${g}'_{ij} (t)$ and add the resulting equations for
$j=1,\cdots ,m$, we find
\begin{equation*}
\begin{split}
 &\frac{1}{2}\partial _t \int_{{\mathbb R}^3} {\theta _{r} (\partial _t
u_{1m} )^2} +\int_{{\mathbb R}^3} {\theta _{r} \partial _t u_{1m}
(\partial _t u_{1m} \partial _{x_1 } u_{1m} +\partial _t u_{2m}
\partial
_{x_2 } u_{1m} +\partial _t u_{3m} \partial _{x_3 } u_{1m} )} + \\
 &\quad \quad +\int_{{\mathbb R}^3} {\theta _{r} \partial _t u_{1m}
(u_{1m} \partial _{x_1 } \partial _t u_{1m} +u_{2m} \partial _{x_2 }
\partial _t u_{1m} +u_{3m} \partial _{x_3 } \partial _t u_{1m} )} + \\
 &\quad \quad +\int_{{\mathbb R}^3} {\theta _{r} \partial _t u_{1m}
\partial _{x_1 } \partial _t p} =\int_{{\mathbb R}^3} {\theta _{r}
\partial _t u_{1m} \,\Delta \partial _t u_{1m} } \\
\end{split}
\end{equation*}
\begin{equation*}
\begin{split}
 &\frac{1}{2}\partial _t \int_{{\mathbb R}^3} {\theta _{r} (\partial _t
u_{2m} )^2} +\int_{{\mathbb R}^3} {\theta _{r} \partial _t u_{2m}
(\partial _t u_{1m} \partial _{x_1 } u_{2m} +\partial _t u_{2m}
\partial
_{x_2 } u_{2m} +\partial _t u_{3m} \partial _{x_3 } u_{2m} )} + \\
 &\quad \quad +\int_{{\mathbb R}^3} {\theta _{r} \partial _t u_{2m}
(u_{1m} \partial _{x_1 } \partial _t u_{2m} +u_{2m} \partial _{x_2 }
\partial _t u_{2m} +u_{3m} \partial _{x_3 } \partial _t u_{2m} )} + \\
 &\quad \quad +\int_{{\mathbb R}^3} {\theta _{r} \partial _t u_{2m}
\partial _{x_2 } \partial _t p} =\int_{{\mathbb R}^3} {\theta _{r}
\partial _t u_{2m} \,\Delta \partial _t u_{2m} } \\
 &\frac{1}{2}\partial _t \int_{{\mathbb R}^3} {\theta _{r} (\partial _t
u_{3m} )^2} +\int_{{\mathbb R}^3} {\theta _{r} \partial _t u_{3m}
(\partial _t u_{1m} \partial _{x_1 } u_{3m} +\partial _t u_{2m}
\partial
_{x_2 } u_{3m} +\partial _t u_{3m} \partial _{x_3 } u_{3m} )} + \\
 &\quad \quad +\int_{{\mathbb R}^3} {\theta _{r} \partial _t u_{3m}
(u_{1m} \partial _{x_1 } \partial _t u_{3m} +u_{2m} \partial _{x_2 }
\partial _t u_{3m} +u_{3m} \partial _{x_3 } \partial _t u_{3m} )} + \\
 &\quad \quad +\int_{{\mathbb R}^3} {\theta _{r} \partial _t u_{3m}
\partial _{x_3 } \partial _t p} =\int_{{\mathbb R}^3} {\theta _{r}
\partial _t u_{3m} \,\Delta \partial _t u_{3m} }  \\
\end{split}
\end{equation*}
\\
\\
Let $ r\to +\infty $,
\begin{equation}
\begin{split}
 &\frac{1}{2}\partial _t \int_{{\mathbb R}^3} {(\partial _t u_{1m} )^2}
+\int_{{\mathbb R}^3} {\partial _t u_{1m} (\partial _t u_{1m}
\partial _{x_1 } u_{1m} +\partial _t u_{2m} \partial _{x_2 } u_{1m}
+\partial _t u_{3m}
\partial _{x_3 } u_{1m} )} + \\
 &\quad \quad +\int_{{\mathbb R}^3} {\partial _t u_{1m} (u_{1m} \partial _{x_1 }
\partial _t u_{1m} +u_{2m} \partial _{x_2 } \partial _t u_{1m} +u_{3m}
\partial _{x_3 } \partial _t u_{1m} )} + \\
 &\quad \quad +\int_{{\mathbb R}^3} {\partial _t u_{1m} \partial _{x_1 } \partial
_t p} =\int_{{\mathbb R}^3} {\partial _t u_{1m} \,\Delta \partial _t u_{1m} } \\
 &\frac{1}{2}\partial _t \int_{{\mathbb R}^3} {(\partial _t u_{2m} )^2}
+\int_{{\mathbb R}^3} {\partial _t u_{2m} (\partial _t u_{1m}
\partial _{x_1 } u_{2m} +\partial _t u_{2m} \partial _{x_2 } u_{2m}
+\partial _t u_{3m}
\partial _{x_3 } u_{2m} )} + \\
 &\quad \quad +\int_{{\mathbb R}^3} {\partial _t u_{2m} (u_{1m} \partial _{x_1 }
\partial _t u_{2m} +u_{2m} \partial _{x_2 } \partial _t u_{2m} +u_{3m}
\partial _{x_3 } \partial _t u_{2m} )} + \\
 &\quad \quad +\int_{{\mathbb R}^3} {\partial _t u_{2m} \partial _{x_2 } \partial
_t p} =\int_{{\mathbb R}^3} {\partial _t u_{2m} \,\Delta \partial _t u_{2m} } \\
 &\frac{1}{2}\partial _t \int_{{\mathbb R}^3} {(\partial _t u_{3m} )^2}
+\int_{{\mathbb R}^3} {\partial _t u_{3m} (\partial _t u_{1m}
\partial _{x_1 } u_{3m} +\partial _t u_{2m} \partial _{x_2 } u_{3m}
+\partial _t u_{3m}
\partial _{x_3 } u_{3m} )} + \\
 &\quad \quad +\int_{{\mathbb R}^3} {\partial _t u_{3m} (u_{1m} \partial _{x_1 }
\partial _t u_{3m} +u_{2m} \partial _{x_2 } \partial _t u_{3m} +u_{3m}
\partial _{x_3 } \partial _t u_{3m} )} + \\
 &\quad \quad +\int_{{\mathbb R}^3} {\partial _t u_{3m} \partial _{x_3 } \partial
_t p} =\int_{{\mathbb R}^3} {\partial _t u_{3m} \,\Delta \partial _t u_{3m} } \\
\end{split}
\end{equation}

Moreover,
\begin{equation*}
\begin{split}
 &\int_{{\mathbb R}^3} {\theta _{r} (\partial _t u_{1m} \partial _{x_1 }
\partial _t p+\partial _t u_{2m} \partial _{x_2 } \partial _t p+\partial _t
u_{3m} \partial _{x_3 } \partial _t p)} \\
 &\quad =-\int_{{\mathbb R}^3} {\theta _{r} \partial _t p\,\;\partial _t
(\partial _{x_1 } u_{1m} +\partial _{x_2 } u_{2m} +\partial _{x_3 }
u_{3m}
)} \\
 &\quad \;\;\;\;  -\int_{{\mathbb R}^3} {\partial _t p\,(\;\partial _t u_{1m}
\partial _{x_1 } \theta _{r} +\partial _t u_{2m} \partial _{x_2 }
\theta _{r} +\partial _t u_{3m} \partial _{x_3 } \theta _{r} )} \\
\end{split}
\end{equation*}
let $ r\to +\infty $ we get
\[
\int_{{\mathbb R}^3} {(\partial _t u_{1m} \partial _{x_1 } \partial
_t p+\partial _t u_{2m} \partial _{x_2 } \partial _t p+\partial _t
u_{3m}
\partial _{x_3 } \partial _t p)} =0
\]
and
\begin{equation*}
\begin{split}
 &\int_{{\mathbb R}^3} {\theta _{r} \partial _t u_{im} (u_{1m} \partial
_{x_1 } \partial _t u_{im} +u_{2m} \partial _{x_2 } \partial _t
u_{im}
+u_{3m} \partial _{x_3 } \partial _t u_{im} )} \\
 &\quad =\frac{1}{2}\int_{{\mathbb R}^3} {\theta _{r} (u_{1m} \partial
_{x_1 } (\partial _t u_{im} )^2+u_{2m} \partial _{x_2 } (\partial _t
u_{im}
)^2+u_{3m} \partial _{x_3 } (\partial _t u_{im} )^2)} \\
 &\quad =-\frac{1}{2}\int_{{\mathbb R}^3} {\theta _{r} (\partial _t
u_{im} )^2(\partial _{x_1 } u_{1m} +\partial _{x_2 } u_{2m}
+\partial _{x_3
} u_{3m} )} \\
 &\quad \;\;\;\;  -\frac{1}{2}\int_{{\mathbb R}^3} {(\partial _t u_{im} )^2(u_{1m}
\partial _{x_1 } \theta _{r} +u_{2m} \partial _{x_2 } \theta
_{r} +u_{3m} \partial _{x_3 } \theta _{r} )} \\
\end{split}
\end{equation*}
let $ r\to +\infty $ we get
\[
\int_{{\mathbb R}^3} {\partial _t u_{im} (u_{1m} \partial _{x_1 }
\partial _t u_{im} +u_{2m} \partial _{x_2 } \partial _t u_{im}
+u_{3m} \partial _{x_3 }
\partial _t u_{im} )} =0,\quad \quad i=1,2,3
\]
\\
as well as
\begin{equation*}
\begin{split}
 &\int_{{\mathbb R}^3} {\theta _{r} \partial _t u_{im} \,\Delta \partial
_t u_{im} } =\int_{{\mathbb R}^3} {\theta _{r} \partial _t u_{im}
(\partial _{x_1 }^2 \partial _t u_{im} +\partial _{x_2 }^2
\partial _t
u_{im} +\partial _{x_3 }^2 \partial _t u_{im} )} \\
 &\quad =-\int_{{\mathbb R}^3} {\theta _{r} ((\partial _{x_1 } \partial
_t u_{im} )^2+(\partial _{x_2 } \partial _t u_{im} )^2+(\partial
_{x_3 }
\partial _t u_{im} )^2)} \\
 &\quad \;\;\;\;   -\int_{{\mathbb R}^3} {\partial _t u_{im} (\,\partial _{x_1 }
\theta _{r} \partial _{x_1 } \partial _t u_{im} +\partial _{x_2 }
\theta _{r} \partial _{x_2 } \partial _t u_{im} +\partial _{x_3 }
\theta _{r} \partial _{x_3 } \partial _t u_{im} )} \\
\end{split}
\end{equation*}
let $ r\to +\infty $ we get
\[
\int_{{\mathbb R}^3} {\partial _t u_{im} \,\Delta \partial _t u_{im}
} =-\int_{{\mathbb R}^3} {((\partial _{x_1 } \partial _t u_{im}
)^2+(\partial _{x_2 } \partial _t u_{im} )^2+(\partial _{x_3 }
\partial _t u_{im} )^2)} ,\quad \quad i=1,2,3
\]
\\
it follows from (19) and above conclusions that
\begin{equation*}
\begin{split}
 &\frac{1}{2}\partial _t \int_{{\mathbb R}^3} {((\partial _t u_{1m} )^2+(\partial
_t u_{2m} )^2+(\partial _t u_{3m} )^2)} \;\; + \\
 &\quad \quad \quad +\left\| {\nabla \partial _t u_{1m} } \right\|_{L^2({\mathbb
R}^3)}^2 +\left\| {\nabla \partial _t u_{2m} }
\right\|_{L^2({\mathbb R}^3)}^2
+\left\| {\nabla \partial _t u_{3m} } \right\|_{L^2({\mathbb R}^3)}^2 \\
 &\quad \le \left\| {\partial _t u_{1m} } \right\|_{L^4({\mathbb R}^3)} \left(
{\left\| {\partial _t u_{1m} } \right\|_{L^4({\mathbb R}^3)} \left\|
{\partial _{x_1 } u_{1m} } \right\|_{L^2({\mathbb R}^3)} +\left\|
{\partial _t u_{2m} } \right\|_{L^4({\mathbb R}^3)} \left\|
{\partial _{x_2 } u_{1m} }
\right\|_{L^2({\mathbb R}^3)} } \right.+ \\
 &\left. {\quad \quad \quad \quad \quad \quad \quad \quad \quad \quad
+\left\| {\partial _t u_{3m} } \right\|_{L^4({\mathbb R}^3)} \left\|
{\partial
_{x_3 } u_{1m} } \right\|_{L^2({\mathbb R}^3)} } \right) \\
 &\quad +\left\| {\partial _t u_{2m} } \right\|_{L^4({\mathbb R}^3)} \left(
{\left\| {\partial _t u_{1m} } \right\|_{L^4({\mathbb R}^3)} \left\|
{\partial _{x_1 } u_{2m} } \right\|_{L^2({\mathbb R}^3)} +\left\|
{\partial _t u_{2m} } \right\|_{L^4({\mathbb R}^3)} \left\|
{\partial _{x_2 } u_{2m} }
\right\|_{L^2({\mathbb R}^3)} +} \right. \\
 &\quad \quad \quad \quad \quad \quad \quad \quad \quad \quad \left.
{+\left\| {\partial _t u_{3m} } \right\|_{L^4({\mathbb R}^3)}
\left\| {\partial
_{x_3 } u_{2m} } \right\|_{L^2({\mathbb R}^3)} } \right) \\
\end{split}
\end{equation*}
\begin{equation*}
\begin{split}
 &\quad +\left\| {\partial _t u_{3m} } \right\|_{L^4({\mathbb R}^3)} \left(
{\left\| {\partial _t u_{1m} } \right\|_{L^4({\mathbb R}^3)} \left\|
{\partial _{x_1 } u_{3m} } \right\|_{L^2({\mathbb R}^3)} +\left\|
{\partial _t u_{2m} } \right\|_{L^4({\mathbb R}^3)} \left\|
{\partial _{x_2 } u_{3m} }
\right\|_{L^2({\mathbb R}^3)} +} \right. \\
 &\quad \quad \quad \quad \quad \quad \quad \quad \quad \quad \left.
{+\left\| {\partial _t u_{3m} } \right\|_{L^4({\mathbb R}^3)}
\left\| {\partial
_{x_3 } u_{3m} } \right\|_{L^2({\mathbb R}^3)} } \right) \\
 &\quad \le \left( {\sum\limits_{i=1}^3 {\left\| {\partial _t u_{im} }
\right\|_{L^4({\mathbb R}^3)}^2 } } \right)^{1/2}\left(
{\sum\limits_{j=1}^3 {\left\| {\partial _t u_{jm} }
\right\|_{L^4({\mathbb R}^3)}^2 } } \right)^{1/2}\left(
{\sum\limits_{i,j=1}^3 {\left\| {\partial _{x_i } u_{jm}
} \right\|_{L^2({\mathbb R}^3)}^2 } } \right)^{1/2} \\
\end{split}
\end{equation*}
\\

Since
\begin{equation*}
\begin{split}
 &\sum\limits_{i=1}^3 {\left\| {\partial _t u_{im} } \right\|_{L^4({\mathbb
R}^3)}^2 } \le 2\sum\limits_{i=1}^3 {\left( {\left\| {\partial _t
u_{im} } \right\|_{L^2({\mathbb R}^3)}^{1/2} \left\| {\nabla
\partial _t u_{im} }
\right\|_{L^2({\mathbb R}^3)}^{3/2} } \right)} \\
 &\quad \le 2\left( {\sum\limits_{i=1}^3 {\left\| {\partial _t u_{im} }
\right\|_{L^2({\mathbb R}^3)}^2 } } \right)^{1/4}\left(
{\sum\limits_{i=1}^3 {\left\| {\nabla \partial _t u_{im} }
\right\|_{L^2({\mathbb R}^3)}^2 } }
\right)^{3/4} \\
\end{split}
\end{equation*}
then
\begin{equation*}
\begin{split}
 &\partial _t \left( {\sum\limits_{i=1}^3 {\left\| {\partial _t u_{im} }
\right\|_{L^2({\mathbb R}^3)}^2 } } \right)+2\left(
{\sum\limits_{i=1}^3 {\left\| {\nabla \partial _t u_{im} }
\right\|_{L^2({\mathbb R}^3)}^2 } }
\right) \\
 &\quad \le 2^2\left( {\sum\limits_{i=1}^3 {\left\| {\partial _t u_{im} }
\right\|_{L^2({\mathbb R}^3)}^2 } } \right)^{1/4}\left(
{\sum\limits_{i=1}^3 {\left\| {\nabla \partial _t u_{im} }
\right\|_{L^2({\mathbb R}^3)}^2 } } \right)^{3/4}\left(
{\sum\limits_{i=1}^3 {\left\| {\nabla u_{im} }
\right\|_{L^2({\mathbb R}^3)}^2 } } \right)^{1/2} \\
 &\quad \le 3^3 \left( {\sum\limits_{i=1}^3 {\left\| {\partial _t u_{im} }
\right\|_{L^2({\mathbb R}^3)}^2 } } \right)\left(
{\sum\limits_{i=1}^3 {\left\| {\nabla u_{im} }
\right\|_{L^2({\mathbb R}^3)}^2 } } \right)^2 + \left(
{\sum\limits_{i=1}^3 {\left\| {\nabla
\partial _t u_{im} }
\right\|_{L^2({\mathbb R}^3)}^2 } } \right) \\
\end{split}
\end{equation*}
\\
it follows that
\[
\partial _t \left( {\sum\limits_{i=1}^3 {\left\| {\partial _t u_{im} }
\right\|_{L^2({\mathbb R}^3)}^2 } } \right) + \left(
{\sum\limits_{i=1}^3 {\left\| {\nabla \partial _t u_{im} }
\right\|_{L^2({\mathbb R}^3)}^2 } } \right) \le \phi _m (t)\left(
{\sum\limits_{i=1}^3 {\left\| {\partial _t u_{im} }
\right\|_{L^2({\mathbb R}^3)}^2 } } \right)
\]
where
\[
\phi _m (t) =  1\,+\, 3^3  \left( {\sum\limits_{i=1}^3 {\left\| {\nabla
u_{im} } \right\|_{L^2({\mathbb R}^3)}^2 } } \right)^2
\]
\\

Introducing a stream function: $\psi =(\psi _2 ,\psi _2 ,\psi _3 )$,
\[
\mbox{curl}\psi =(\partial _{x_2 } \psi _3 -\partial _{x_3 } \psi _2
,\;\,\;\partial _{x_3 } \psi _1 -\partial _{x_1 } \psi _3 ,\;\,\;\partial
_{x_1 } \psi _2 -\partial _{x_2 } \psi _1 )
\]
According to $\omega =\mbox{curl}u$, $u=\mbox{curl}\psi $ and
$\mbox{div}\psi =0$, we have
\[
\mbox{curlcurl}\psi =-\Delta \psi =\omega ,
\quad
-\Delta \mbox{curl}\psi =\mbox{curl}\omega
\]
That is, $-\Delta u=\mbox{curl}\omega $. Then $(-\Delta
u,\;\,u)=(\mbox{curl}\omega ,\;\,u)$, where
\[
(-\Delta u,\;\,\theta _{r} u)=\sum\limits_{i=1}^3 {(-\Delta u_i
,\;\,\theta _{r} u_i )} =\sum\limits_{i=1}^3 {(\nabla u_i
,\;\,\theta _{r} \nabla u_i )} +\sum\limits_{i=1}^3 {(\nabla u_i
,\;\,u_i \nabla \theta _{r} )}
\]
let $ r\to +\infty $ we get
\[
(-\Delta u,\;\,u)=\sum\limits_{i=1}^3 {(\nabla u_i ,\;\,\nabla u_i
)} =\sum\limits_{i=1}^3 {\left\| {\nabla u_i }
\right\|_{L^2({\mathbb R}^3)}^2 }
\]
\\

In addition,
\begin{equation*}
\begin{split}
 &(\mbox{curl}\omega ,\;\,\theta _{r} u)=(\partial _{x_2 } \omega _3
-\partial _{x_3 } \omega _2 ,\;\;\theta _{r} u_1 )+(\partial _{x_3
} \omega _1 -\partial _{x_1 } \omega _3 ,\;\;\theta _{r} u_2 ) \\
 &\quad \quad \quad \quad \quad \quad \quad +(\partial _{x_1 } \omega _2
-\partial _{x_2 } \omega _1 ,\;\;\theta _{r} u_3 ) \\
 &\quad =-(\omega _3 ,\;\theta _{r} \partial _{x_2 } u_1 )+(\omega
_2 ,\;\theta _{r} \partial _{x_3 } u_1 )-(\omega _1 ,\;\theta
_{r} \partial _{x_3 } u_2 ) \\
 &\qquad  +(\omega _3 ,\;\theta _{r} \partial _{x_1 } u_2
)-(\omega _2 ,\;\theta _{r} \partial _{x_1 } u_3 )+(\omega _1
,\;\theta _{r} \partial _{x_2 } u_3 ) \\
 &\qquad  -(\omega _3 ,\;u_1 \partial _{x_2 } \theta _{r}
)+(\omega _2 ,\;u_1 \partial _{x_3 } \theta _{r} )-(\omega _1
,\;u_2 \partial _{x_3 } \theta _{r} ) \\
 &\qquad  +(\omega _3 ,\;u_2 \partial _{x_1 } \theta _{r}
)-(\omega _2 ,\;u_3 \partial _{x_1 } \theta _{r} )+(\omega _1
,\;u_3 \partial _{x_2 } \theta _{r} ) \\
\end{split}
\end{equation*}
let $ r\to +\infty $ we get
\begin{equation*}
\begin{split}
 &(\mbox{curl}\omega ,\;\,u)=-(\omega _3 ,\;\partial _{x_2 } u_1 )+(\omega _2
,\;\partial _{x_3 } u_1 )-(\omega _1 ,\;\partial _{x_3 } u_2 )+(\omega _3
,\;\partial _{x_1 } u_2 ) \\
 &\qquad \quad \quad \quad \quad -(\omega _2 ,\;\partial _{x_1 } u_3 )+(\omega
_1 ,\;\partial _{x_2 } u_3 ) \\
 &\quad =(\omega _1 ,\;\;\partial _{x_2 } u_3 -\partial _{x_3 } u_2 )+(\omega
_2 ,\;\;\partial _{x_3 } u_1 -\partial _{x_1 } u_3 )+(\omega _3
,\;\;\partial _{x_1 } u_2 -\partial _{x_2 } u_1 ) \\
 &\quad =(\omega ,\;\,\mbox{curl}u)=(\omega ,\omega )=\sum\limits_{i=1}^3
{\left\| {\omega _i } \right\|_{L^2({\mathbb R}^3)}^2 } \\
\end{split}
\end{equation*}
Hence
\[
\left( {\sum\limits_{i=1}^3 {\left\| {\nabla u_i }
\right\|_{L^2({\mathbb R}^3)}^2 } } \right)^{1/2}=\left(
{\sum\limits_{i=1}^3 {\left\| {\omega _i } \right\|_{L^2({\mathbb
R}^3)}^2 } } \right)^{1/2}
\]
it follows that
\[
\phi _m (t) = 1\,+\, 3^3 \left( {\sum\limits_{i=1}^3 {\left\| {\omega _{im} }
\right\|_{L^2({\mathbb R}^3)}^2 } } \right)^2<+\infty
\]
\\

By the Gronwall inequality,
\[
\frac{d}{dt}\left\{ {\left( {\sum\limits_{i=1}^3 {\left\| {\partial
_t u_{im} } \right\|_{L^2({\mathbb R}^3)}^2 } } \right)\;\exp \left(
{-\int_0^t {\phi _m (s)ds} } \right)} \right\}\le 0
\]
whence
\[
\mathop {\sup }\limits_{t\in (0,T)} \left( {\sum\limits_{i=1}^3
{\left\| {\partial _t u_{im} (t)} \right\|_{L^2({\mathbb R}^3)}^2 }
} \right)\le \left( {\sum\limits_{i=1}^3 {\left\| {\partial _t
u_{im} (0)} \right\|_{L^2({\mathbb R}^3)}^2 } } \right)\;\exp \left(
{\int_0^T {\phi _m (s)ds} } \right)
\]
Therefore
\[
\partial _t u_{im} \in L^\infty (0,T;\;H)\cap L^\infty (0,T;\;V),\quad \quad
i=1,2,3
\]
\\

Finally we write (1) in the form
\[
\sum\limits_{i=1}^3 {(-\Delta (\,\theta _{r} u_i ),\;v_i )}
=\sum\limits_{i=1}^3 {(-\theta _{r} \partial _t u_i -\theta _{r}
(u\cdot \nabla )u_i +g_i ,\;\;v_i )} ,\quad \quad v_i \in V
\]
where
\[
g_i = -\; u_i \,\Delta \theta _{r} -\;2\,(\nabla \theta _{r}
,\;\,\nabla u_i ) + p\,\partial _{x_i } \theta _{r}
\]
That is,
\[
\sum\limits_{i=1}^3 {(\,\nabla (\,\theta _{r} u_i ),\;\;\nabla v_i
)} =\sum\limits_{i=1}^3 {(-\theta _{r} \partial _t u_i -\theta _{r}
(u\cdot \nabla )u_i +g_i ,\;\;v_i )}
\]

Since
\[
\partial _t u_i \in L^\infty (0,T;\;H),\quad \quad (u\cdot \nabla )u_i \in
L^\infty (0,T;\;H)
\]
Similar to the Theorem 3.8 in Chapter 3 of [4], and let $ r\to
+\infty $, we obtain
\[
u_i \in L^\infty (0,T;\;H^2({\mathbb R}^3)),\quad \quad i=1,2,3
\]
\\
\\

\textbf{Remark 1. }Noting that $(-\Delta u,\;v)=(-\partial _t
u-(u\cdot \nabla )u,\;v)$. Since $\partial _t u$ and $(u\cdot \nabla
)u$ are of some degree of continuity, then $u$ can reach a higher
degree of continuity, based upon the smoothing effect of inverse
elliptic operator $\Delta ^{-1}$. By repeated application of this
process one can prove that the solution $u$ is in $C^\infty
((0,T)\times {\mathbb R}^3)$.
\\

\textbf{Remark 2. } Based on problems separated and potential theory
of fluid flow, we may keep the same result for the general
initial-boundary value problems of 3D Navier-Stokes equation under
the assumptions of regularity on the boundary and data.
\\
\\
\\
\\

\end{document}